\documentclass{amsart}
\usepackage{amsmath,mathtools}
\usepackage{amsrefs}
\usepackage{microtype}
\usepackage{lmodern}
\usepackage{hyperref}
\usepackage{tikz-cd}
\usepackage{graphicx}
\usepackage{tikz-network}
\usetikzlibrary{arrows}
\usepackage{enumerate}
 
\usetikzlibrary{decorations.markings}
\usepackage{caption}
\theoremstyle{plain}
\newtheorem{theorem}[subsection]{Theorem}
\newtheorem{proposition}[subsection]{Proposition}
\newtheorem{lemma}[subsection]{Lemma}

\theoremstyle{definition}
\newtheorem{definition}[subsection]{Definition}

\theoremstyle{remark}
\newtheorem{remark}[subsection]{Remark}

\numberwithin{equation}{section}

\newcommand{\cal}{\mathcal}
\newcommand{\cla}{{\cal A}}

\newcommand{\cle}{{\cal E}}
\newcommand{\clf}{{\cal F}}

\newcommand{\der}{\textup{d}}

\def\a*{{\cal A}_{h,*}}
\def\B{{\cal B}(h)}
\def\B1{{\cal B}_1(h)}
\def\b{{\cal B}^{\rm s.a.}(h)}
\def\b1{{\cal B}^{\rm s.a.}_1(h)}

\newcommand{\ot}{ \ \otimes \ }

\newcommand{\raro}{\rightarrow}

\newcommand{\id}{\mbox{id}}

\def\a*{{\cal A}_{h,*}}
\def\B1{{\cal B}_1(h)}
\def\b{{\cal B}^{\rm s.a.}(h)}
\def\b1{{\cal B}^{\rm s.a.}_1(h)}

\usepackage{graphicx} 
\begin{document}

\title[]{Almost complex structure on finite points from bidirected graphs}
\author[Joardar]{Soumalya Joardar}
\address{Department of Mathematics and Statistics, Indian Institute of Science Education and Research Kolkata, Mohanpur - 741246, West
	Bengal, India}
\email{soumalya@iiserkol.ac.in}
\author[Rahaman]{Atibur Rahaman}
\address{Department of Mathematics and Statistics, Indian Institute of Science Education and Research Kolkata, Mohanpur - 741246, West
	Bengal, India}
\email{atibur.pdf@iiserkol.ac.in}
\begin{abstract}
We show that there is an almost complex structure on a differential calculus on finite points coming from a {\bf bidirected} finite graph without multiple edges or loops. We concentrate on a polygon as a concrete case. In particular, a `holomorphic structure on the exterior bundle' built from the polygon is studied. Also a positive Hochschild $2$-cocycle on the vertex set of the polygon, albeit a trivial one, is shown to arise naturally from the almost complex structure.  
\end{abstract}
\maketitle
\section{Introduction}
Despite much progress in noncommutative geometry over the past few decades, noncommutative complex geometry is still believed to be in a relatively early phase of its development. Unlike the classical counterpart where there is a nice harmony between the analytic and algebraic aspects of complex geometry, the notions in noncommutative complex geometry are less settled. In \cite{Connescuntz}, a possible approach to complex structure on noncommutative geometry has been speculated in terms of a positive Hochschild cocycle on an involutive algebra. There are much more `algebraic approach' to noncommutative complex geometry in the literature (for example see \cites{Beggs, qbundle, Heckenbergkolb, Schwartz, Staford}). In \cite{NCtori} it is shown that the almost complex structure on noncommutative torus in \cite{Schwartz} defines a positive Hochschild $2$-cocycle, whereas in \cite{qbundle}, the almost complex structure gives rise to a twisted positive Hochschild $2$-cocycle.\\
\indent As noncommutative complex geometry is mostly in its early developmental phase, awaiting analogues to classical theorems like the Kodaira embedding theorem (see section 5.3 of \cite{Huybrechts}) or Serre's GAGA, it becomes imperative to construct new examples. These examples would play a crucial role in testing any potential approach toward noncommutative complex geometry.  An almost complex structure on the quantum group analogues of irreducible flag manifolds has been studied in \cites{Heckenbergkolb, Heckenbergkolb2, Heckenbergkolb3}. In the past few years the study of noncommutative complex geometry has gained some momentum. Readers should see \cites{kahlerbuachala,Das} for most recent developments about K$\ddot{{\rm a}}$hler structures on some noncommutative complex manifolds. But so far the examples are mostly restricted to quantum analogues of flag manifolds or noncommutative tori. More recently, K$\ddot{{\rm a}}$hler structures have been discovered in the setting of holomorphic etále groupoid in \cite{debashishda}. In this context, from a complementary viewpoint regarding examples, we explore an almost complex structure on finite points. To be more precise, we study almost complex structure on finitely many points arising from bidirected finite graphs. We take an algebraic approach and in particular we follow the notion of an almost complex structure from \cite{Beggs} (or equivalently \cites{qbundle, Heckenbergkolb}). The $\ast$-differential calculi on finite points is well studied. For example in \cite{Majid}, it is shown that the category of first order differential calculus is the same as the category of directed finite graphs (without loops or multiple edges) on finitely many points. Furthermore, such differential calculi coming from directed graphs are inner and can be extended via a `canonical prolongation procedure' to obtain forms of all degrees. We take this as our starting point and show in this paper that under suitable conditions, the differential calculus can be given an almost complex structure provided the graph is {\bf bidirected}. The idea is to define a $\mathbb{C}$-linear map $J$ on the first order differential calculus using the bidirectedness and subsequently extend it to the spaces of all forms. The integrability of such an almost complex structure follows from the inner structure of the differential calculus coming from bidirected graphs. Once integrable almost complex structure is obtained, one can study the relevant complex geometric notions. We concentrate on a concrete model where the $\ast$-differential calculus comes from a bidirected $n$-gon and the resulting calculus is of dimension $2$. The spaces of one and two forms are shown to be free modules so that one can consider `holomorphic sections' of the `exterior bundle' with respect to a natural choice of holomorphic structure. The holomorphic sections form a ring under wedge product and is shown to be isomorphic to the exterior algebra $\Lambda^{\bullet}(\mathbb{C}^{2})$. Moreover, we construct a positive Hochschild $2$-cocycle on the underlying vertex set of the polygon out of the almost complex structure following the footsteps of noncommutative $2$-tori in \cite{NCtori}.  \vspace{0.05in}\\
{\bf Notations and conventions}: We use $\bar{\ot}$ to denote the interior tensor product of two bimodules. For a bimodule $\cle$, the $n$-fold interior tensor product of $\cle$ will be denoted by $\cle^{\bar{\otimes} n}$. ${\rm Sp}_{\mathbb{C}}$ stands for complex linear span. 

\section{Almost complex structure}
\subsection{Generalities}

\begin{definition}(\cite{Khalkhali})
Let $\cla$ be a $\ast$-algebra over the field of complex numbers $\mathbb{C}$.
Then a $\ast$-differential calculus is a data $(\Omega^{\bullet}(\cla),\der,\ast)$  where
\begin{enumerate}[(i)]
    \item $\Omega^{\bullet}(\cla)=\oplus_{n\geq 0}\Omega^{n}(\cla)$ is a graded $\ast$-algebra such that $\Omega^{\bullet}(\cla)$ is also an $\cla$\ndash$\cla$-bimodule with $\Omega^{0}(\cla)=\cla$ by convention;
    \item $\der:\Omega^{\bullet}(\cla)\raro\Omega^{\bullet+1}(\cla)$ is a graded derivation i.e. $\der$ is a $\mathbb{C}$-linear map satisfying the graded Leibniz rule
    \begin{displaymath}
    \der(\omega\wedge\eta)=\der\omega\wedge\eta+(-1)^{{\rm deg}(\omega)}\omega\wedge \der\eta,
    \end{displaymath}
    for $\omega,\eta\in \Omega^{\bullet}(\cla)$ and $\der^2=0$;
    \item For all $\omega\in\Omega^{\bullet}(\cla)$, $\der(\omega^{\ast})=(\der\omega)^{\ast}$.
\end{enumerate}
\end{definition} 

\begin{remark}
 If $\omega\in\Omega^{n}(\cla)$, the degree of $\omega$ denoted above by ${\rm deg}(\omega)$ is equal to $n$. We have denoted the product of the graded algebra $\Omega^{\bullet}(\cla)$ by $\wedge$. With these notations, $\ast$ is an antilinear graded involution meaning,
\begin{displaymath}
(\omega\wedge\eta)^{\ast}=(-1)^{{\rm deg}(\omega){\rm deg}(\eta)}\eta^{\ast}\wedge\omega^{\ast}.
\end{displaymath}
Also if ${\rm deg}(\omega)=k$, ${\rm deg}(\eta)=l$, then ${\rm deg}(\omega\wedge\eta)=k+l$.
\end{remark}
\begin{definition}
A $\ast$-differential calculus over a $\ast$-algebra $\cla$ is said to be of dimension $k$ if $\Omega^{k}(\cla)\neq 0$, and $\Omega^{n}(\cla)=0$ for all $n\geq k+1$.     
\end{definition}

\begin{definition}
    A $\ast$-differential calculus of dimension $k$ is said to be orientable if there is an $\cla$\ndash$\cla$ bimodule isomorphism between the bimodule of top forms $\Omega^{k}(\cla)$ and $\cla$. A choice of such an isomorphism is called an orientation.
\end{definition}
\begin{definition}( \cite{Khalkhali}*{Def. 3.6.3})
    \label{defgradedtrace}
     A closed graded trace of dimension $k$ on a differential calculus $(\Omega^{\bullet}(\cla),\der,\ast)$ is a $\mathbb{C}$-linear map $\int:\Omega^{k}(\cla)\raro\mathbb{C}$ such that
     \begin{enumerate}[(i)]
        \item $\int \der\omega=0$ for all $\omega\in\Omega^{k-1}$;
        \item For $\omega_{1}\in\Omega^{i}(\cla), \omega_{2}\in\Omega^{j}(\cla)$ such that $i+j=k$,
            \begin{displaymath}
              \int(\omega_{1}\wedge\omega_{2}-(-1)^{ij}\omega_{2}\wedge\omega_{1})=0.
            \end{displaymath}
    \end{enumerate}
\end{definition}

In this paper we shall construct a $\ast$-differential calculus over a $\ast$-algebra $\cla$ from an `inner first order' $\ast$-differential calculus $(\Omega^{1}(\cla),\der,\ast)$ by a canonical prolongation procedure as in \cite{CMQG}. Recall that for an inner first order $\ast$-differential calculus $(\Omega^{1}(\cla),\der,\ast)$, $\Omega^{1}(\cla)$ is an $\cla-\cla$-bimodule such that $\der:\cla\raro\Omega^{1}(\cla)$ is a $\mathbb{C}$-linear derivation i.e. $\der(ab)=a\der{b}+\der{a}b$ for all $a,b\in\cla$, and $\der(\omega)^{\ast}=\der(\omega^{\ast})$ for all $\omega\in\Omega^{1}(\cla)$. The calculus is inner means that there is an element $\theta\in\Omega^{1}(\cla)$ such that $\der{a}=[\theta,a]$ for all $a\in\cla$ where $[\theta,a]$ is defined by the formula $[\theta,a]:=\theta.a-a.\theta$ using the bimodule structure of $\Omega^{1}(\cla)$. Let us discuss the canonical prolongation procedure as described in \cite{CMQG}. To that end let $(\Omega^{1}(\cla),\der,\ast)$ be an inner first order $\ast$-differential calculus over a $\ast$-algebra $\cla$ and $\sigma:\Omega^{1}(\cla)\bar{\ot} \Omega^{1}(\cla)\raro\Omega^{1}(\cla)\bar{\ot}\Omega^{1}(\cla)$ be a bimodule map such that $\sigma$ satisfies the following braid relation on $\Omega^{1}(\cla)^{\bar{\otimes}3}$:
\begin{displaymath}
\sigma_{12}\sigma_{23}\sigma_{12}=\sigma_{23}\sigma_{12}\sigma_{23},
\end{displaymath}
where $\sigma_{23}=({\rm id}\bar{\ot}\sigma), \sigma_{12}=(\sigma\bar{\ot}{\rm id})$. Let $n$ be a natural number with $n\geq 2$. Any permutation $p$ on $n$ symbols can be written as product of nearest neighbor transpositions. Let $p=t_{k_{1}}\cdots t_{k_{I(p)}}$ where $I(p)$ is the number of inversions in the sequence $p(1), p(2), \ldots , p(n)$ and $t_{i}$'s are nearest neighbor transpositions as in \cite{CMQG}*{Eq. (3.20)}. For $p$, $\Pi_{p}$ is defined as $\sigma_{k_{1}}\sigma_{k_{2}}\cdots \sigma_{k_{I(p)}}$. The bimodule maps $\sigma_{k}$'s are defined for $k=1,2,\ldots,n-1$ as in~\cite{CMQG}*{Eq. (3.17)}. Then the antisymmetrization map $A_{n}$ is defined as
\begin{displaymath}
A_{n}=\sum_{p\in S_{n}} {\rm sgn}(p)\Pi_{p}.
\end{displaymath}
 Clearly $A_{n}$ is a bimodule map as each $\Pi_{p}$ is so. Then the bimodule of $n$-forms $\Omega^{n}(\cla)$ is defined to be $\Omega^{1}(\cla)^{\bar{\otimes}n}/{\rm Ker}(A_{n})$. We shall denote the map from $\Omega^{1}(\cla)^{\bar{\otimes}n}$ to $\Omega^{n}(\cla)$ mapping an element $\omega\in\Omega^{1}(\cla)^{\bar{\otimes}n}$ to its class $[\omega]\in\Omega^{n}(\cla)$ by $\wedge$. Observe that the antisymmetrization map does not depend on a particular representation of a permutation thanks to the braid relation of $\sigma$. Note that from~\cite{CMQG}*{Eq. (3.25)}, it can be deduced that if $\Omega^{n}(\cla)=0$ for some $n=k$, then $\Omega^{n}(\cla)=0$ for all $n\geq k$. We record the following lemma whose proof is obvious. 
\begin{lemma}
\label{star_2form}
Let $(\Omega^{1}(\cla),\der,\ast)$ be an inner first order $\ast$-differential calculus over a $\ast$-algebra $\cla$. Then there is a well defined antilinear graded involution $\ast$ on $\Omega^{1}(\cla)\bar{\ot}\Omega^{1}(\cla)$ given by $(\omega_{1}\bar{\ot}\omega_{2})^{\ast}=-\omega_{2}^{\ast}\bar{\ot}\omega_{1}^{\ast}$.
\end{lemma}
We prove the following theorem which says that in fact the $\ast$ can be extended to all forms under a suitable condition.  Note that this result differs slightly from~\cite{CMQG}*{Theorem 3.4} and the inner structure is not really needed to prove the result. Although we remark that for our application purpose, we consider the bimodule map $\sigma$ to be such that $\sigma^{2}={\rm id}$ and therefore the proof of the following theorem reduces to that of \cite{CMQG}*{Theorem 3.4} due to~\cite{CMQG}*{Eq. (3.37)}. We decided to keep the general result nonetheless.  
\begin{theorem}
\label{astallforms}
Let $\Omega^{\bullet}(\cla):=\oplus_{n\geq 0}\Omega^{n}(\cla)$ be the bimodule of forms of all degree obtained as a canonical prolongation with respect to a bimodule map $\sigma$. If  $\ast$ on $\Omega^{1}(\cla)\bar{\ot}\Omega^{1}(\cla)$ as defined in Lemma \ref{star_2form} commutes with the bimodule map $\sigma$, then the involution $\ast$ on $\Omega^{1}(\cla)$ extends as a graded antilinear involution on $\Omega^{\bullet}(\cla)$.
\end{theorem}
\begin{proof}
Let $n\geq 2$. Following \cite{CMQG}, for $\eta=\omega_{1}\bar{\ot}\cdots\bar{\ot}\omega_{n}\in\Omega^{1}(\cla)^{\bar{\otimes}n}$, we define
\begin{displaymath}
    \eta^{\ast}:=(-1)^{\frac{n(n-1)}{2}}\omega_{n}^{\ast}\bar{\ot}\cdots\bar{\ot}\omega_{1}^{\ast}.
\end{displaymath}
As observed in\cite{CMQG}, this is clearly an antilinear involution on $\Omega^{1}(\cla)^{\bar{\otimes}n}$. To show that the above definition descends to $\Omega^{n}(\cla)$, we have to show that for $\eta\in{\rm Ker}(A_{n})$, $\eta^{\ast}\in{\rm Ker}(A_{n})$. In fact we shall show that $(A_{n}(\eta^{\ast}))=(A_{n}(\eta))^{\ast}$, proving the claim. Therefore, we take an element $\eta=\omega_{1}\bar{\ot}\cdots\bar{\ot}\omega_{n}$ and observe that for $j<n$,
\begin{eqnarray*}
  \begin{aligned}
    &\sigma_{j(j+1)}(\eta^{\ast})\\
    &\quad= (-1)^{\frac{n(n-1)}{2}}\sigma_{j(j+1)}(\omega_{n}^{\ast}\bar{\ot}\cdots\bar{\ot}\omega_{n-j+1}^{\ast}\bar{\ot}\omega_{n-j}^{\ast}\bar{\ot}\cdots\bar{\ot}\omega_{1}^{\ast})\\
    &\quad= (-1)^{\frac{n(n-1)}{2}}\omega_{n}^{\ast}\bar{\ot}\cdots\bar{\ot}\sigma(\omega_{n-j+1}^{\ast}\bar{\ot}\omega_{n-j}^{\ast})\bar{\ot}\cdots\bar{\ot}\omega_{1}^{\ast},
  \end{aligned}
\end{eqnarray*}
which is equal to 
\begin{displaymath}
   (-1)^{\frac{n(n-1)}{2}+1}\omega_{n}^{\ast}\bar{\ot}\cdots\bar{\ot}(\sigma(\omega_{n-j}\bar{\ot}\omega_{n-j+1}))^{\ast}\bar{\ot}\cdots\bar{\ot}\omega_{1}^{\ast}, 
\end{displaymath}
since, by assumption, $\sigma$ commutes with $\ast$ on $\Omega^{1}(\cla)\bar{\ot}\Omega^{1}(\cla)$. But by definition of $\ast$, the last expression is equal to 
\begin{displaymath}
    (\sigma_{(n-j)(n-j+1)}(\omega_{1}\bar{\ot}\cdots\bar{\ot}\omega_{n}))^{\ast}.
\end{displaymath}
Therefore, for any $j<n$, $\sigma_{j(j+1)}(\eta^{\ast})=(\sigma_{(n-j)(n-j+1)}(\eta))^{\ast}$. Now take any transposition $p=(k \ l)\in S_{n}$ with $k<l$. Then $(k \ l)$ can be written as a product of neighboring transpositions in the following way:
\begin{displaymath}
    (k \ l)=(k \ k+1)(k+1 \ k+2)\cdots(l-1 \ l)(l-2 \ l-1)\cdots(k \ k+1).
\end{displaymath}
Hence by definition of $\Pi_{p}$,
\begin{eqnarray*}
    \Pi_{p}(\eta^{\ast})&=&\sigma_{k(k+1)}\sigma_{(k+1)(k+2)}\cdots\sigma_{(l-1)l}\sigma_{(l-2)(l-1)}\cdots\sigma_{k(k+1)}(\eta^{\ast})\\
    &=&\sigma_{k(k+1)}\sigma_{(k+1)(k+2)}\cdots\sigma_{(l-1)l}\sigma_{(l-2)(l-1)}\cdots\bigl(\sigma_{(n-k)(n-k+1)}(\eta)\bigr)^{\ast}\\
    &=& \biggl(\sigma_{(n-k)(n-k+1)}\sigma_{(n-k-1)(n-k)}\cdots\sigma_{(n-l+1)(n-l)}\cdots\sigma_{(n-k)(n-k+1)}(\eta)\biggr)^{\ast}.
\end{eqnarray*}
But $\sigma_{(n-k)(n-k+1)}\sigma_{(n-k-1)(n-k)}\cdots\sigma_{(n-l+1)(n-l)}\cdots\sigma_{(n-k)(n-k+1)}=\Pi_{p^{\prime}}$, where $p^{\prime}$ is the transposition $(n-k+1 \ n-l+1)$. Now with similar computation one can show that for any cycle $p=(i_{1} \ i_{2} \ \cdots \ i_{k})$,
\begin{displaymath}
    \Pi_{p}(\eta^{\ast})=(\Pi_{p^{\prime}}(\eta))^{\ast},
\end{displaymath}
where $p^{\prime}$ is the cycle $(n-i_{1}+1 \ \cdots \ n-i_{k}+1)$. Writing down an arbitrary permutation as a product of disjoint cycles, the association $p\rightleftharpoons p^{\prime}$ between the cycles clearly extends to an one-one onto correspondence between all permutations so that $\Pi_{p}(\eta^{\ast})=(\Pi_{p^{\prime}}(\eta))^{\ast}$. Moreover, clearly ${\rm sgn}(p)={\rm sgn}(p^{\prime})$. Therefore for any $\eta\in \Omega^{1}(\cla)^{\bar{\otimes}n}$,
\begin{eqnarray*}
    A_{n}(\eta^{\ast})&=&\sum_{p\in S_{n}}{\rm sgn}(p)\Pi_{p}(\eta^{\ast})\\
    &=& \biggl(\sum_{p^{\prime}\in S_{n}}{\rm sgn}(p^{\prime})\Pi_{p^{\prime}}(\eta)\biggr)^{\ast}\\
    &=& \bigl(A_{n}(\eta)\bigr)^{\ast},
\end{eqnarray*}
proving the theorem.
\end{proof}
 The following theorem completes the construction of the canonical prolongation procedure whose proof is a straightforward adaptation of \cite{CMQG}*{Theorem 4.1} and hence omitted.
\begin{theorem}
\label{diffcalculus}
Let $(\Omega^{1}(\cla),\der,\ast)$ be an inner first order $\ast$-differential calculus over a $\ast$-algebra $\cla$ with $\der{a}=[\theta,a]$ for all $a\in\cla$ for some $\theta\in\Omega^{1}(\cla)$. Let $\sigma:\Omega^{1}(\cla)\bar{\ot}\Omega^{1}(\cla)\raro \Omega^{1}(\cla)\bar{\ot}\Omega^{1}(\cla)$ be a bimodule map obeying the braid relation such that $\sigma(\theta\bar{\ot}\theta)=\theta\bar{\ot}\theta$ and $\sigma$ commutes with $\ast$ on $\Omega^{1}(\cla)\bar{\ot}\Omega^{1}(\cla)$. Then there is a unique $\mathbb{C}$-linear graded derivation $\der\colon\Omega^{\bullet}(\cla)\rightarrow\Omega^{\bullet+1}$ extending the derivation $\der:\cla\rightarrow\Omega^{1}(\cla)$. The map $\der$ is given by the following formula for $\omega\in\Omega^{n}(\cla)$:
\begin{equation}
    \label{dformula}
    \der{\omega}:=[\theta,\omega]_{{\rm grad}}= \theta\wedge\omega-(-1)^{n}\omega\wedge\theta.
\end{equation}
Moreover, $(\der\omega)^{\ast}=\der(\omega^{\ast})$ for $\omega\in\Omega^{\bullet}(\cla)$ making $(\Omega^{\bullet}(\cla),\der,\ast)$ a $\ast$-differential calculus over $\cla$.
\end{theorem}
We shall now recall some generalities of almost complex structure on a $\ast$-differential calculus over a $\ast$-algebra. These general facts regarding an almost complex structure are mostly taken from \cite{Beggs}.
\begin{definition}
Let $(\Omega^{\bullet}(\cla),\der,\ast)$ be a $\ast$-differential calculus over a $\ast$-algebra $\cla$. An almost complex structure on $(\Omega^{\bullet}(\cla),\der,\ast)$ is a degree zero derivation $J:\Omega^{\bullet}(\cla)\rightarrow\Omega^{\bullet}(\cla)$ such that
\begin{enumerate}[1.]
    \item $J$ is identically zero on $\cla$ and hence an $\cla$\ndash$\cla$-bimodule endomorphism of $\Omega^{\bullet}(\cla)$;
    \item $J^{2}=-{\rm Id}$ on $\Omega^{1}(\cla)$; and
    \item $J(\omega^{\ast})=(J\omega)^{\ast}$ for all $\omega\in\Omega^{1}(\cla)$.
\end{enumerate}
\end{definition}
Given an almost complex structure on a $\ast$-differential calculus, the bimodule of one-forms has the following decomposition into eigenspaces of $J$:\begin{displaymath}
\Omega^{1}(\cla)=\Omega^{1,0}(\cla)\oplus\Omega^{0,1}(\cla),
\end{displaymath}
where $\Omega^{1,0}(\cla):=\{\omega\in\Omega^{1}(\cla): J(\omega)=i\omega\}$, $\Omega^{0,1}:=\{\omega\in\Omega^{1}(\cla): J(\omega)=-i\omega\}$. $\Omega^{1,0}(\cla)$ and $\Omega^{0,1}(\cla)$ are again $\cla$\ndash$\cla$-bimodules and called the bimodules of $(1,0)$ and $(0,1)$-forms respectively. More generally, one has the higher dimensional analogues of the $(1,0)$ and $(0,1)$-forms. For a fixed $n$, the bimodule $\Omega^{n}(\cla)$ has the following decomposition into bimodules of $(p,q)$-forms for $p+q=n$, where $p,q$ are non negative integers(\cite{Beggs}*{Sec. 2.5}): 
\begin{displaymath}
\Omega^{n}(\cla)=\oplus_{p+q=n}\Omega^{p,q}(\cla),
\end{displaymath}
where $\Omega^{p,q}(\cla):=\{\omega\in\Omega^{n}(\cla):J(\omega)=i(p-q)\omega\}$. Moreover, 
\begin{eqnarray}
\label{integrability}
\Omega^{p,q}(\cla)\wedge\Omega^{p^{\prime},q^{\prime}}(\cla)\subset\Omega^{p+p^{\prime},q+q^{\prime}}(\cla).
\end{eqnarray}
Using the decomposition of $\Omega^{n}(\cla)$ into $(p,q)$-forms, one can define $\partial:\Omega^{p,q}(\cla)\raro\Omega^{p+1,q}$ and $\overline{\partial}:\Omega^{p,q}(\cla)\raro\Omega^{p,q+1}$ in the following way:
\indent For all non-negative integers $p,q$, consider the projections $\pi^{p,q}:\Omega^{p+q}(\cla)\raro\Omega^{p,q}(\cla)$ associated to the decomposition of $\Omega^{p+q}(\cla)$. Then $\partial:\Omega^{p,q}(\cla)\raro\Omega^{p+1,q}(\cla)$ and $\overline{\partial}:\Omega^{p,q}(\cla)\raro\Omega^{p,q+1}$ are defined as 
\begin{displaymath}
\partial:=\pi^{p+1,q}\circ \der, \qquad \overline{\partial}:=\pi^{p,q+1}\circ \der.
\end{displaymath}

\begin{definition} (see~\cite{Beggs}*{Definition 3.1})
An almost complex structure on a $\ast$-differential calculus is said to be integrable if any of the equivalent conditions of the following Lemma \ref{intstructure} holds.
\end{definition}

\begin{lemma} \textup{(}\cite{Beggs}*{Lemma 3.2}\textup{)}
\label{intstructure} Let $J$ be an almost complex structure on a $\ast$-differential calculus $(\Omega^{\bullet}(\cla),\der,\ast)$. Then the following conditions are equivalent:
\begin{enumerate}[1.]
    \item $\partial^{2}:\cla\raro\Omega^{2}(\cla)$ is $0$ as an operator.
    \item $\overline{\partial}^{2}:\cla\raro\Omega^{2}(\cla)$ is $0$ as an operator.
    \item $\der(\Omega^{1,0}(\cla))\subset \Omega^{2,0}(\cla)\oplus\Omega^{1,1}(\cla)$.
    \item $\der(\Omega^{0,1}(\cla))\subset \Omega^{1,1}(\cla)\oplus\Omega^{0,2}(\cla)$.
\end{enumerate}
\end{lemma}
Let $J$ be an integrable almost complex structure on a $\ast$-differential calculus $(\Omega^{\bullet}(\cla),\der,\ast)$. Then thanks to~\cite{Beggs}*{Proposition 3.6}, one can consider the following $p$-th Dolbeault complex for each non negative integer $p$:
\begin{align}
\begin{array}{ccccccccc}
  0 &\longrightarrow & \Omega^{p,0}(\cla) & \overset{\overline{\partial}}{\longrightarrow} & \Omega^{p,1}(\cla) & \overset{\overline{\partial}}{\longrightarrow} & \Omega^{p,2}(\cla) & \overset{\overline{\partial}}\longrightarrow & ...
  \end{array}
  \end{align}
  The $q$-th cohomology of the above sequence is called the $(p,q)$ Dolbeault cohomology of $\cla$ and is denoted by $H_{\overline{\partial}}^{p,q}(\cla)$. 
  \subsection{Holomorphic modules}
  \begin{definition} (\cite{Beggs})
      Let $(\Omega^{\bullet}(\cla),\der,\ast)$ be a $\ast$-differential calculus with an integrable almost complex structure $J$. A left $\overline{\partial}$-operator on a left $\cla$-module $\cle$ is a $\mathbb{C}$-linear map $\overline{\nabla}:\cle\raro\Omega^{0,1}(\cla)\bar{\ot}\cle$ such that
      \begin{displaymath}
          \overline{\nabla}(a.e)=\overline{\partial}a\bar{\ot}e+a\overline{\nabla}(e),
      \end{displaymath}
      for $a\in\cla, e\in\cle$.
  \end{definition}
  For $q\geq 1$, one can define $\overline{\nabla}:\Omega^{0,q}(\cla)\bar{\ot}\cle\raro\Omega^{0,q+1}\bar{\ot}\cle$ by
  \begin{displaymath}
  \overline{\nabla}(\omega\bar{\ot}e):=\overline{\partial}\omega\bar{\ot} e+(-1)^{q}\omega\wedge\overline{\nabla}(e).
  \end{displaymath}
  The holomorphic curvature is defined to be $\overline{\nabla}^{2}:\cle\raro\Omega^{0,2}(\cla)\bar{\ot}\cle$.
  \begin{definition}
      Let $(\Omega^{\bullet}(\cla),\der,\ast)$ be a $\ast$-differential calculus with an integrable almost complex structure $J$. A holomorphic structure on a left module $\cle$ is a $\overline{\partial}$-operator $\overline{\nabla}$ on $\cle$ such that the holomorphic curvature $\overline{\nabla}^{2}$ vanishes.
  \end{definition}
  With a holomorphic structure on a left module $\cle$, $(\cle,\overline{\nabla})$ is called a holomorphic left module. Given a holomorphic left module, one can consider the following complex:
  \begin{align}
\begin{array}{ccccccccc}
  0 &\longrightarrow & \cle & \overset{\overline{\nabla}}{\longrightarrow} & \Omega^{0,1}(\cla)\bar{\ot}\cle & \overset{\overline{\nabla}}{\longrightarrow} & \Omega^{0,2}(\cla)\bar{\ot}\cle & \overset{\overline{\nabla}}\longrightarrow & ...
  \end{array}
  \end{align}
  The elements of the $0$-th cohomology group $H^{0}(\cle,\overline{\nabla})$ of the complex are the holomorphic sections of the bundle represented by the module $\cle$.  
\begin{remark}
  To call $\cle$ the module of sections of a `noncommutative vector bundle' one generally adds the technical demand that $\cle$ should be finitely generated and  projective. In our case, this technical condition will be satisfied as we shall study the holomorphic sections of the modules which will be free.  
\end{remark}
    \indent Recall that a left connection on a left $\cla$-module $\cle$ is a $\mathbb{C}$-linear map $\nabla:\cle\raro\Omega^{1}(\cla)\bar{\ot}\cle$ such that 
  \begin{displaymath}
  \nabla(a.e)=\der{a}\bar{\ot}e+a.\nabla(e).
  \end{displaymath}
  Given a left connection the curvature is $\nabla^{2}:\cle\raro\Omega^{2}(\cla)\bar{\ot}\cle$. We can induce a holomorphic structure on a left module $\cle$ from a left connection whose $\Omega^{0,2}$ curvature component vanishes. More precisely, we state the following without proof:
  \begin{proposition} (\cite{Beggs})
  Let $(\Omega^{\bullet}(\cla),\der,\ast)$ be a $\ast$-differential calculus with an integrable almost complex structure and $\cle$ be a left $\cla$-module having a left connection $\nabla$. Then the map
  \begin{displaymath}
      \overline{\nabla}:=(\pi^{0,1}\ot {\rm id})\circ\nabla:\cle\raro\Omega^{0,1}(\cla)\bar{\ot}\cle
  \end{displaymath}
  is a left $\overline{\partial}$-operator. Moreover, if the $\Omega^{0,2}$-component of $\nabla^2$ vanishes, then $(\cle,\overline{\nabla})$ is a holomorphic module. 
  \end{proposition}
 Now we recall the notion of bimodule connections. Again as before in the following $(\Omega^{\bullet},\der,\ast)$ is a $\ast$-differential calculus.
 \begin{definition} (\cite{qbundle})
     \label{bimoduleconnection}
A bimodule connection on an $\cla$-bimodule $\cle$ with respect to a bimodule isomorphism $\psi:\cle\bar{\ot}\Omega^{1}(\cla)\raro\Omega^{1}(\cla)\bar{\ot}\cle$ is a left connection $\nabla$ on $\cle$ with such that following twisted Leibniz rule is satisfied:   
\begin{displaymath}
    \nabla(e.a)=\nabla(e).a+\psi(e\bar{\ot}\der{a}).
\end{displaymath}
 \end{definition}
 As noted in \cite{qbundle}, the above definition applies to the differential calculus $(\Omega^{0,1}(\cla),\overline{\partial})$ giving a notion of holomorphic bimodules. Given bimodule connections $(\nabla_{1},\nabla_{2})$ on bimodules $\cle,\clf$ with respect to bimodule isomorphisms $\psi_{1},\psi_{2}$ respectively, one can consider their tensor product bimodule connection $\nabla$ on the bimodule $\cle\bar{\ot}\clf$ with respect to the bimodule isomorphism $(\psi_{1}\bar{\ot}{\rm id})\circ({\rm id}\bar{\ot}\psi_{2})$ given by (see \cite{qbundle}):
 \begin{displaymath}
     \nabla:=(\nabla_{1}\bar{\ot} 1)+(\psi_{1}\bar{\ot} 1)(1\bar{\ot}\nabla_{2}).
 \end{displaymath}
  \indent Returning to our initial framework where we have an inner first order $\ast$-differential calculus $(\Omega^{1}(\cla),\der,\ast)$ over a $\ast$-algebra $\cla$, with the introduction of some additional data, one can systematically construct a canonical prolongation of the first order calculus to obtain a $\ast$-differential calculus over $\cla$. If we have a map $J$ on $\Omega^{1}(\cla)$, then under specific conditions governing $J$, it becomes possible to equip the $\ast$-differential calculus with an integrable complex structure which effectively `extends' the structure on $\Omega^{1}(\cla)$. To see this, first let us assume that there is an $\cla$\ndash$\cla$-bimodule map $J:\Omega^{1}(\cla)\raro\Omega^{1}(\cla)$. Then we record the following lemma whose proof is obvious.
  \begin{lemma}
  \label{LemmaJextension2form}
  The map $J$ extends to a map from $\Omega^{1}(\cla)\bar{\ot}\Omega^{1}(\cla)$ to $\Omega^{1}(\cla)\bar{\ot}\Omega^{1}(\cla)$ given by the formula
  \begin{eqnarray}  \label{Jextension}  J(\omega_{1}\bar{\ot}\omega_{2})=J(\omega_{1})\bar{\ot}\omega_{2}+\omega_{1}\bar{\ot} J(\omega_{2}).
  \end{eqnarray}
      \end{lemma}
      Now we are ready to state and prove the theorem which gives us a condition governing $J$ so that $J$ extends to a canonical prolongation to give an integrable complex structure.
  \begin{theorem}
  \label{extension}
  Let $(\Omega^{\bullet}(\cla),\der,\ast)$ be a $\ast$-differential calculus obtained as a canonical prolongation of an {\bf inner} first order $\ast$-differential calculus $(\Omega^{1}(\cla),\der,\ast)$ corresponding to a bimodule map $\sigma:\Omega^{1}(\cla)\bar{\ot}\Omega^{1}(\cla)\raro\Omega^{1}(\cla)\bar{\ot}\Omega^{1}(\cla)$. Assume that there is an $\cla$\ndash$\cla$-bimodule map $J:\Omega^{1}(\cla)\raro\Omega^{1}(\cla)$ satisfying the following:
  \begin{enumerate}[1.]
     \item $J^2=-{\rm Id}$; $J(\omega^{\ast})=(J\omega)^{\ast}$.
     \item The extension $J:\Omega^{1}(\cla)\bar{\ot}\Omega^{1}(\cla)\raro\Omega^{1}(\cla)\bar{\ot}\Omega^{1}(\cla)$ given by Lemma \ref{LemmaJextension2form}, commutes with the bimodule map $\sigma$.
    \end{enumerate}
   Then $J$ extends to a bimodule endomorphism on $\Omega^{\bullet}(\cla)$ which is also a degree zero derivation inducing an {\bf integrable} almost complex structure on $(\Omega^{\bullet}(\cla),\der,\ast)$. 
  \end{theorem}
  
  \begin{proof}
      For any $n$, we define $J:\Omega^{1}(\cla)^{\bar{\otimes}n}\raro\Omega^{1}(\cla)^{\bar{\otimes}n}$ by
      \begin{displaymath}
          J(\omega_{1}\bar{\ot}\omega_{2}\bar{\ot}\cdots\bar{\ot}\omega_{n}):=\sum_{k=1}^{n}\omega_{1}\bar{\ot}\omega_{2}\bar{\ot}\cdots\bar{\ot} J(\omega_{k})\bar{\ot}\cdots\bar{\ot}\omega_{n}.          
      \end{displaymath}
      Clearly $J$ is a well defined degree zero derivation on $\Omega^{1}(\cla)^{\bar{\otimes}n}$. $J$ will descend to a degree zero derivation on $\Omega^{n}(\cla)$ if $J$ preserves ${\rm ker}(A_{n})$. We claim that for any neighboring transposition $(k \ k+1)\in S_{n}$, $J\biggl(\sigma_{k(k+1)}(\omega_{1}\bar{\ot}\cdots\bar{\ot}\omega_{n})\biggr)=\sigma_{k(k+1)}\biggl(J(\omega_{1}\bar{\ot}\cdots\bar{\ot}\omega_{n})\biggr)$. We shall prove the claim for the transposition $(1\ 2)$. The claim for the other transpositions can be proved similarly. 
      \begin{eqnarray*}
       \begin{aligned}
        &J\biggl(\sigma_{12}(\omega_{1}\bar{\ot}\cdots\bar{\ot}\omega_{n})\biggr)\\
        &\quad= J\biggl(\sigma(\omega_{1}\bar{\ot}\omega_{2})\bar{\ot}\omega_{3}\bar{\ot}\cdots\bar{\ot}\omega_{n}\biggr)\\
        &\quad= J\bigl(\sigma(\omega_{1}\bar{\ot}\omega_{2})\bigr)\bar{\ot}\omega_{3}\bar{\ot}\cdots\bar{\ot}\omega_{n}\\
        &\qquad+\sum_{k=3}^{n}\sigma(\omega_{1}\bar{\ot}\omega_{2})\bar{\ot}\omega_{3}\bar{\ot}\cdots\bar{\ot}J(\omega_{k})\bar{\ot}\cdots\bar{\ot}\omega_{n}\\
        &\quad= \sigma_{12}\biggl(J(\omega_{1}\bar{\ot}\omega_{2}\bar{\ot}\cdots\bar{\ot}\omega_{n})\biggr).
        \end{aligned}
      \end{eqnarray*}
      In the last equality we have used the fact that $J$ commutes with $\sigma$ on $\Omega^{1}(\cla)\bar{\ot}\Omega^{1}(\cla)$. Therefore, from the definition of the antisymmterization map, it is clear that $A_{n}(J\eta)=J(A_{n}\eta)$ for any $\eta\in\Omega^{1}(\cla)^{\bar{\otimes}n}$ proving that $J$ descends to a well defined degree zero derivation on $\Omega^{n}(\cla)$ for any $n$ inducing an almost complex structure on the $\ast$-differential calculus $(\Omega^{\bullet}(\cla),\der,\ast)$. To show that the structure is integrable, recall that the calculus is a canonical prolongation of an inner first order differential calculus and therefore for any $\omega\in\Omega^{1,0}(\cla)$, $\der\omega=\theta\wedge\omega+\omega\wedge\theta$ for some $\theta\in\Omega^{1}(\cla)$. Hence by Equation \eqref{integrability}, $\der\omega\subset\Omega^{2,0}(\cla)\oplus\Omega^{1,1}(\cla)$. So by the equivalent condition (3) of Lemma \ref{intstructure}, the almost complex structure is integrable.
  \end{proof}
  \begin{remark}
  The integrability condition is satisfied automatically since we are considering an inner differential calculus.
  \end{remark}
  \subsection{Bidirected graphs}
  Given a finite set $V$, it is well known (see \cites{Majid,Connes}) that there is a one to one correspondence between the differential structures of $C(V)$ and directed graphs (without multiple edges or loops) with $V$ as the set of vertices. For us, a graph will always mean a graph {\bf without multiple edges or loops}. Given a directed graph over $V$, if one denotes the edge set by $E$, then the vector space of functions on $E$ becomes a $C(V)$\ndash$C(V)$-bimodule with the following bimdule structure:\\
  \indent Let $x,y\in V$. The notation $x\raro y\in E$ will mean that there is a directed edge from $x$ to $y$, with source being $x$ (written as $s(e)=x$) and target being $y$ (written as $t(e)=y$). The function on $E$ which takes the value $1$ on $x\raro y$ and $0$ on the rest of the edges will be denoted by $\xi_{x\raro y}$. Then clearly $\{\xi_{x\raro y}\}_{x\raro y\in E}$ forms a $\mathbb{C}$-linear basis of $C(E)$. With these notations, the left and right action of $C(V)$ on the basis elements of $C(E)$ are given by
  \begin{eqnarray}
  f.\xi_{x\raro y}=f(x)\xi_{x\raro y}, \ \xi_{x\raro y}.f=\xi_{x\raro y}f(y),
  \end{eqnarray}
  for $f\in C(V)$. There is a derivation $\der:C(V)\raro C(E)$ given by (\cite{Majid})
  \begin{eqnarray}
  \der{f}:=\sum_{x\raro y\in E}(f(y)-f(x))\xi_{x\raro y}.
  \end{eqnarray}
  Thus $(C(E),C(V), \der)$ becomes a first order differential calculus on $C(V)$. From now on we shall denote $C(V)$ by $\cla$ and $C(E)$ by $\Omega^{1}(\cla)$. Note that $\cla$ is clearly a $\ast$-algebra with pointwise complex conjugation. As observed in \cite{Majid}, this first order differential calculus is inner with $\der{f}=[\theta,f]$ where $\theta=\sum_{x\raro y\in E}\xi_{x\raro y}$.\\ 
  \indent To discuss about an almost complex structure we shall need our graphs to be bidirected. So from now on we shall only consider {\bf bidirected graphs without multiple edges or loops}. Therefore let $(E,V)$ be a finite bidirected graph over the vertex set $V$. Then one has the first order differential calculus as discussed earlier. We can define an antilinear involution on $\Omega^{1}(\cla)$ making $(\Omega^{1}(\cla),\der)$ a first order $\ast$-calculus. The antilinear involution on the linear basis is given by
  \begin{displaymath}
  (\xi_{x\raro y})^{\ast}:=-\xi_{y\raro x}.
  \end{displaymath}
  Observe that bidirectedness is used to define the involution. One can easily see that $(\Omega^{1}(\cla),\der,\ast)$ is a first order $\ast$-differential calculus over the $\ast$-algebra $\cla=C(V)$. Now let us discuss a canonical prolongation of the calculus $(\Omega^{1}(\cla),\der,\ast)$. To begin, note that the $\cla$\ndash$\cla$-bimodule $\Omega^{1}(\cla)\bar{\ot}\Omega^{1}(\cla)$ possesses a $\mathbb{C}$- linear basis $\{\xi_{x\raro y}\bar{\ot}\xi_{y\raro z}\}_{x\raro y, y\raro z\in E}$. This follows from the obvious formula (see \cite{Majid}):
  \begin{displaymath}
  \xi_{x\raro y}\bar{\ot}\xi_{z\raro w}=\delta_{y,z}\xi_{x\raro y}\bar{\ot}\xi_{y\raro w}.
  \end{displaymath}
  Now to perform the canonical prolongation procedure to obtain a $\ast$-differential calculus $(\Omega^{\bullet}(\cla),\der,\ast)$, by Theorem \ref{diffcalculus}, one needs a bimodule map $\sigma:\Omega^{1}(\cla)\bar{\ot}\Omega^{1}(\cla)\raro \Omega^{1}(\cla)\bar{\ot}\Omega^{1}(\cla)$ such that\\ (i) $\sigma$ satisfies the braid relation on $\Omega^{1}(\cla)^{\bar{\otimes}3}$.\\
  (ii) $\sigma(\theta\bar{\ot}\theta)=(\theta\bar{\ot}\theta)$, where $\theta=\sum_{x\raro y\in E}\xi_{x\raro y}$.\\
  (iii) $\sigma$ commutes with the $\ast$-map on $\Omega^{1}(\cla)\bar{\ot}\Omega^{1}(\cla)$ as mentioned in Lemma \ref{star_2form}. \\
  \indent As $\{\xi_{x\raro y}\bar{\ot}\xi_{y\raro z}\}_{x\raro y, y\raro z\in E}$ is a linear basis for $\Omega^{1}(\cla)\bar{\ot}\Omega^{1}(\cla)$, any bimodule map $\sigma$ can be written as (see \cite{Majid})
  \begin{displaymath}
  \sigma(\xi_{x\raro y}\bar{\ot}\xi_{y\raro z})=\sum_{w}\sigma^{x,y,z}_{w}\xi_{x\raro w}\bar{\ot}\xi_{w\raro z},
  \end{displaymath}
  for some scalars $\sigma^{x,y,z}_{w}\in\mathbb{C}$. Guided by \cite{Majid}, we say a bimodule map is of permutation type if for fixed $x,z\in V$ such that there is at least one two-length path from $x$ to $z$, the matrix $\sigma^{x,\bullet,z}_{\bullet}$ is a permutation matrix. It is an easy observation (see \cite{Majid}) that a bimodule map $\sigma$ of permutation type satisfies $\sigma(\theta\bar{\ot}\theta)=(\theta\bar{\ot}\theta)$. In order to fulfill condition (iii) above i.e. to make $\sigma$ commute with $\ast$, we shall restrict our attention to the following special class of bimodule maps:
  \begin{definition}
  Let $(V,E)$ be a finite, bidirected graph. Then we call a bimodule map $\sigma$ of permutation type invariant under path reversal if for $x,z\in V$ such that $x\neq z$, the  permutation matrix $\sigma^{x,\bullet,z}_{\bullet}$ is the same as the permutation matrix $\sigma^{z,\bullet,x}_{\bullet}$.
  \end{definition}
  \begin{remark}
  Note that the assumption of bidirectedness is crucial for the above definition to make sense. The bidirectedness ensures that for $x\neq z$, there is a two-length path between $x,z$ iff there is a two length path between $z,x$ obtained by reversing the arrows of the paths.
  \end{remark}
  In the following $(V,E)$ is a finite, bidirected, graph so that it has an inner first order $\ast$-differential calculus $(\Omega^{1}(\cla),\der,\ast)$.
  \begin{theorem}
      \label{prolongationgraph}
      Let $\sigma:\Omega^{1}(\cla)\bar{\ot}\Omega^{1}(\cla)\raro\Omega^{1}(\cla)\bar{\ot}\Omega^{1}(\cla)$ be a bimodule map of permutation type which is invariant under path reversal and satisfies the braid relation on $\Omega^{1}(\cla)^{\bar{\otimes}3}$. Then $\sigma$ commutes with the antilinear involution $\ast$ as defined in Lemma \ref{star_2form}. Consequently, there is a canonical prolongation $(\Omega^{\bullet}(\cla),\der,\ast)$ of $(\Omega^{1}(\cla),\der,\ast)$ with respect to $\sigma$.
  \end{theorem}
  \begin{proof}
  The calculus $(\Omega^{1}(\cla),\der,\ast)$ is inner with $\der{f}=[\theta,f]$ where $\theta=\sum_{x\raro y\in E}\xi_{x\raro y}$. Since $\sigma$ is of permutation type, $\sigma(\theta\bar{\ot}\theta)=(\theta\bar{\ot}\theta)$. By assumption $\sigma$ satisfies the braid relation. Therefore the proof will be complete by Theorem \ref{diffcalculus} if $\sigma$ commutes with $\ast$ on $\Omega^{1}(\cla)\bar{\ot}\Omega^{1}(\cla)$. As observed earlier, $\Omega^{1}(\cla)\bar{\ot}\Omega^{1}(\cla)$ has a $\mathbb{C}$-linear basis $\{\xi_{x\raro y}\bar{\ot}\xi_{y\raro z}\}_{x\raro y, y\raro z\in E}$. The graph being bidirected, for each $x\in V$ there are always basis elements of the form $\{\xi_{x\raro y}\bar{\ot}\xi_{y\raro x}\}_{x\raro y\in E}$. On such basis elements, by definition of $\ast$,
  \begin{eqnarray*}
  \sigma(\xi_{x\raro y}\bar{\ot}\xi_{y\raro x})^{\ast}=-\sigma(\xi_{x\raro y}\bar{\ot}\xi_{y\raro x})=-\xi_{x\raro z}\bar{\ot}\xi_{z\raro x},
  \end{eqnarray*}
  for some $z$ such that $x\raro z\in E$ which is easily seen to be equal to $(\sigma(\xi_{x\raro y}\bar{\ot}\xi_{y\raro x}))^{\ast}$. For the basis elements of the form $\xi_{x\raro y}\bar{\ot}\xi_{y\raro z}$ where $x\neq z$,
  \begin{eqnarray*}
  \sigma(\xi_{x\raro y}\bar{\ot}\xi_{y\raro z})^{\ast}=-\sigma(\xi_{z\raro y}\bar{\ot}\xi_{y\raro x})=-\xi_{z\raro y^{\prime}}\bar{\ot}\xi_{y^{\prime}\raro x}.
  \end{eqnarray*}
  On the other hand, as $\sigma$ is invariant under path reversal,
  \begin{eqnarray*}
  (\sigma(\xi_{x\raro y}\bar{\ot}\xi_{y\raro z}))^{\ast}=(\xi_{x\raro y^{\prime}}\bar{\ot}\xi_{y^{\prime}\raro z})^{\ast}=-\xi_{z\raro y^{\prime}}\bar{\ot}\xi_{y^{\prime}\raro x}.
  \end{eqnarray*}
  Therefore using $\mathbb{C}$-linearity of $\sigma$, $(\sigma(\omega))^{\ast}=\sigma(\omega^{\ast})$ for all $\omega\in\Omega^{1}(\cla)\bar{\ot}\Omega^{1}(\cla)$ completing the proof of the theorem.
  \end{proof}
  From now on, given a bidirected, finite graph $(V,E)$ we shall consider the above $\ast$-differential calculus obtained as a canonical prolongation of the first order $\ast$-differential calculus $(\Omega^{1}(\cla),\der,\ast)$. This calculus is associated with the graph structure with respect to a bimodule map \(\sigma\) of permutation type, subject to the conditions of being invariant under path reversal and satisfying the braid relation.
  Note that choice of such $\sigma$ is far from unique. We shall focus on the problem of defining an integrable almost complex structure on such a $\ast$-differential calculus. To achieve this, let's begin with the first order $\ast$-differential calculus $(\Omega^{1}(\cla),\der,\ast)$ given by the graph structure. Our objective is to define a bimodule map $J:\Omega^{1}(\cla)\raro\Omega^{1}(\cla)$ satisfying $J^{2}=-{\rm Id}$ and $J(\omega^{\ast})=(J(\omega))^{\ast}$. To proceed, consider two adjacent vertices $x,y\in V$. Given the bidirected nature of the graph, there are two edges $x\raro y$ and $y\raro x$. Define $J$ by 
  \begin{eqnarray}
   &&J(\xi_{x\raro y})=i\xi_{x\raro y}\\
   && J(\xi_{y\raro x})=-i\xi_{y\raro x}.
  \end{eqnarray} Considering all the adjacent vertices we define $J$ on all the $\mathbb{C}$-linear basis elements of $\Omega^{1}(\cla)$ and extend $\mathbb{C}$-linearly to all of $\Omega^{1}(\cla)$. It is clear from the definition that $J^{2}=-{\rm Id}$ and $J(\omega^{\ast})=(J\omega)^{\ast}$. 
 \begin{remark}
  Again like the bimodule map $\sigma$, the choice of $J$ is far from unique. However, the condition $J(\omega^{\ast})=(J\omega)^{\ast}$ puts the restriction that if $\xi_{x\raro y}$ is in the $i$-th eigen space, then $\xi_{y\raro x}$ must be in the $-i$-th eigen space of $J$. The point is that one always has at least one choice of $J$ on bidirected graphs.   
 \end{remark}  This map $J$ extends to a derivation from $\Omega^{1}(\cla)\bar{\ot}\Omega^{1}(\cla)$ to $\Omega^{1}(\cla)\bar{\ot}\Omega^{1}(\cla)$ by Lemma \ref{LemmaJextension2form}. Therefore if $J$ commutes with $\sigma$ then the differential calculus $(\Omega^{\bullet}(\cla),\der,\ast)$ being inner, $J$ induces an {\bf integrable} almost complex structure by Theorem \ref{extension}. We shall call such an integrable almost complex structure {\bf an integrable almost complex structure coming from the graph $(V,E)$}. We shall not pursue the general conditions on $\sigma$ or $J$ so that $J$ commutes with $\sigma$. Instead we shall consider the  bidirected polygon on $n$-points as a concrete case and discuss an integrable almost complex structure on it. We shall do it in the next section. At this point, it is worth mentioning that the condition that $J$ commutes with $\sigma$, puts some further restriction on the choice of $J$.\\
 \indent Before going to the next section, let us briefly discuss holomorphic structure on the trivial $\cla$\ndash$\cla$ bimodule $\cla(=C(V))$ with respect to an integrable almost complex structure $(\Omega^{\bullet}(\cla),\der,\ast,J)$ coming from a bidirected graph $(V,E)$. Recall that such a structure in particular produces a first order $\ast$-differential calculus $(\Omega^{0,1}(\cla),\overline{\partial},\ast)$ and hence a directed graph over $V$. We call it the holomorphic part of the original bidirected graph. $\overline{\partial}$ is itself a canonical bimodule $\overline{\partial}$-operator on the bimodule $\cla$ with respect to the isomorphism $\psi(\equiv {\rm id}):\Omega^{0,1}(\cla)\cong\cla\bar{\ot}\Omega^{0,1}\raro\Omega^{0,1}(\cla)\bar{\ot}\cla\cong \Omega^{0,1}(\cla)$. As $\overline{\partial}^2=0$, $(\cla,\overline{\partial})$ is a holomorphic bimodule. 
  \begin{lemma}
      \label{uniqueholomorphic}
      For a bidirected graph $(V,E)$, let $(\Omega^{\bullet}(\cla),\der,\ast,J)$ be an integrable almost complex structure coming from $(V,E)$. Then $\overline{\partial}$ gives a unique bimodule holomorphic structure to the trivial bimodule $\cla$. Moreover, if the holomorphic part is connected, then $H^{0}(\cla,\overline{\partial})=\mathbb{C}$.
  \end{lemma}
  \begin{proof}
      Let $\nabla$ be a bimodule $\overline{\partial}$-operator on $\cla$. Then $\Phi\equiv(\nabla-\overline{\partial}):\cla\raro\Omega^{0,1}(\cla)\bar{\ot}\cla\cong\Omega^{0,1}(\cla)$ is an $\cla$\ndash$\cla$ bimodule map. Therefore for any $x\in V$, $\Phi(\delta_{x})=\sum_{k,l}\lambda_{k,l}\xi_{k\raro l}$, for some scalars $\lambda_{k,l}$ where ${k\raro l}$ are the edges of the holomorphic part of the bidirected graph $(V,E)$. Then for $y\neq x$,
      \begin{displaymath}          
      \Phi(\delta_{x}\delta_{y})=\Phi(\delta_{x})\delta_{y}=\Phi(\delta_{y}\delta_{x})=\delta_{y}\Phi(\delta_{x})=0.
      \end{displaymath}
      But $\Phi(\delta_{x})\delta_{y}=\sum_{k:k\raro y}\lambda_{k,y}\xi_{k\raro y}$. Therefore for any $k$ such that there is an edge from $k$ to $y$, $\lambda_{k,y}=0$. Similarly using left $\cla$-linearity of $\Phi$, for any $k$ such that there is an edge from $y$ to $k$, $\lambda_{y,k}=0$. Varying over all $y\neq x$ in $V$, we conclude that $\lambda_{i,j}=0$ if $i$ or $j$ is different from $x$. As the graph doesn't have any loop, $\Phi(\delta_{x})=0$. As $x$ was arbitrary, $\Phi\equiv 0$, proving that $\nabla=\overline{\partial}$. The fact that $H^{0}(\cla,\overline{\partial})=\mathbb{C}$ can be proved along the lines of \cite{Majid}. $H^{0}(\cla,\overline{\partial})$ is the $\mathbb{C}$-linear span of the constant function ${\bf 1}$.
  \end{proof}
  
\section {A case study: Polygon on finite points}

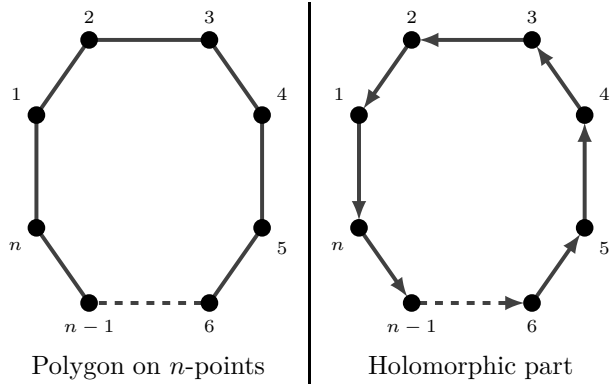
\begin{figure}[h]
    \begin{minipage}[h]{0.5\textwidth}
    \begin{tabular}{@{}c|c@{}}
    \begin{tikzpicture}
 \Vertex[label=$1$,position=above left,shape=circle,size=0.05,color=black]{1}
  \Vertex[x=.7, y=1, label=$2$,position=above,shape=circle, size=0.05,color=black]{2}
  \Vertex[x=2.3,y=1,label=$3$, position=above,shape=circle, size=0.05,color=black]{3}
  \Vertex[x=3,label=$4$,position=above right,shape=circle, size=0.05,color=black]{4}
  \Vertex[,y=-1.5,label=$n$,position=below left,shape=circle, size=0.05,color=black]{8}
  \Vertex[x=3,y=-1.5,label=$5$,position=below right,shape=circle, size=0.05,color=black]{5}
  \Vertex[x=2.3,y=-2.5,position=below,label=$6$,shape=circle, size=0.05,color=black]{6}
  \Vertex[x=.7,y=-2.5,label=$n-1$,position=below,shape=circle, size=0.05,color=black]{7}

  \Edge[](1)(2);
  \Edge[](2)(3);
  \Edge[](3)(4);
  \Edge[](4)(5);
  \Edge[](5)(6);
  \Edge[style={dashed}](6)(7);
  \Edge[](1)(8);
  \Edge[](7)(8);
\end{tikzpicture}
&
\begin{tikzpicture}
\Vertex[label=$1$,position=above left,shape=circle,size=0.05,color=black]{1}
  \Vertex[x=.7, y=1, label=$2$,position=above,shape=circle, size=0.05,color=black]{2}
  \Vertex[x=2.3,y=1,label=$3$, position=above,shape=circle, size=0.05,color=black]{3}
  \Vertex[x=3,label=$4$,position=above right,shape=circle, size=0.05,color=black]{4}
  \Vertex[,y=-1.5,label=$n$,position=below left,shape=circle, size=0.05,color=black]{8}
  \Vertex[x=3,y=-1.5,label=$5$,position=below right,shape=circle, size=0.05,color=black]{5}
  \Vertex[x=2.3,y=-2.5,position=below,label=$6$,shape=circle, size=0.05,color=black]{6}
  \Vertex[x=.7,y=-2.5,label=$n-1$,position=below,shape=circle, size=0.05,color=black]{7}
  \Edge[Direct](2)(1);
  \Edge[Direct](3)(2);
  \Edge[Direct](4)(3);
  \Edge[Direct](5)(4);
  \Edge[Direct](6)(5);
  \Edge[style={dashed},Direct](7)(6);
  \Edge[Direct](1)(8);
  \Edge[Direct](8)(7);
\end{tikzpicture}
\\
Polygon on \(n\)-points & Holomorphic part
    \end{tabular}
\end{minipage}
\end{figure}

Let $(V,E)$ be an $n$-gon. We label the vertices by integers $1,2,\ldots,n$. We consider it as a bidirected graph. In the following we shall construct a $\ast$-differential calculus over the $\ast$-algebra $C(V)$ and denote $C(V)$ by $\cla$. The graph $(V,E)$ being  bidirected, induces a first order $\ast$-differential calculus $(\Omega^{1}(\cla),\der,\ast)$. To construct the canonical prolongation we need a bimodule map $\sigma$ with suitable properties. Before constructing a canonical prolongation let us clarify notations to be used throughout the rest of the paper. Whenever an arrow $\mu\raro \nu$ represents a directed edge in the polygon, $\mu,\nu$ are chosen from the additive group $\mathbb{Z}_{n}=\{1,2,\ldots,n\}$ with $n$ being the additive identity of $\mathbb{Z}_{n}$. For example for $\mu=1$, $\mu\raro\mu-1$ will represent the directed edge from vertex $1$ to vertex $n$. Likewise for $\mu=n$, $\mu\raro\mu+1$ will represent the directed edge from vertex $n$ to vertex $1$. With this notation, we start with the observation that $\Omega^{1}(\cla)$ has the following basis:
\begin{equation*}
    \{\xi_{\mu\raro\mu+1}, \xi_{\mu\raro\mu-1}\}_{\mu=1,2,\ldots,n}.
\end{equation*} 
Now we shall construct the bimodule map $\sigma$ to perform a canonical prolongation of the first order $\ast$-differential calculus. Note that $\Omega^{1}(\cla)\bar{\ot}\Omega^{1}(\cla)$ has a $\mathbb{C}$-linear basis consisting of elements of the following form for $\mu=1,2,\ldots,n$:
\begin{eqnarray}\label{omega2basis}
&&\xi_{\mu\raro{\mu+1}}\bar{\ot}\xi_{\mu+1\raro\mu+2}, \nonumber \\
&&\xi_{\mu\raro{\mu-1}}\bar{\ot}\xi_{\mu-1\raro\mu-2}, \nonumber \\
&&\xi_{\mu\raro\mu+1}\bar{\ot}\xi_{\mu+1\raro\mu},\\
&&\xi_{\mu\raro\mu-1}\bar{\ot}\xi_{\mu-1\raro\mu}. \nonumber 
\end{eqnarray}
We define $\sigma$ on the above basis elements by the following formulae and extend it $\mathbb{C}$-linearly to the whole of $\Omega^{1}(\cla)\bar{\ot}\Omega^{1}(\cla)$:
\begin{eqnarray*}
&&\sigma(\xi_{\mu\raro\mu+1}\bar{\ot}\xi_{\mu+1\raro\mu})=\xi_{\mu\raro\mu-1}\bar{\ot}\xi_{\mu-1\raro\mu},\\
&&\sigma(\xi_{\mu\raro\mu-1}\bar{\ot}\xi_{\mu-1\raro\mu})=\xi_{\mu\raro\mu+1}\bar{\ot}\xi_{\mu+1\raro\mu},
\end{eqnarray*}
and 
$\sigma$ acts trivially on the rest of the basis elements. 
\begin{lemma}
    \label{sigmaok}
     $\sigma$ is a bimodule map of permutation type invariant under path reversal. Moreover, $\sigma$ satisfies the braid relation on $\Omega^{1}(\cla)^{\bar{\otimes}3}$.
\end{lemma}
\begin{proof}
The fact that $\sigma$ is a bimodule map of permutation type invariant under path reversal is easy to check. To check the braid relation, it is enough to check the relation on the vector space basis elements of $\Omega^{1}(\cla)^{\bar{\otimes}3}$. A basis for $\Omega^{1}(\cla)^{\bar{\otimes}3}$ is given by collection of three length paths emitting from each vertex of the $n$-gon. One can list them down for $\mu=1,\ldots,n$:
    \begin{eqnarray}  \label{omega3basis}  &&\xi_{\mu\raro\mu+1}\bar{\ot}\xi_{\mu+1\raro\mu+2}\bar{\ot}\xi_{\mu+2\raro\mu+3},\nonumber \\
    && \xi_{\mu\raro\mu-1}\bar{\ot}\xi_{\mu-1\raro\mu-2}\bar{\ot}\xi_{\mu-2\raro\mu-3},\nonumber \\
    && \xi_{\mu\raro\mu+1}\bar{\ot}\xi_{\mu+1\raro\mu}\bar{\ot}\xi_{\mu\raro\mu+1},\nonumber \\
    && \xi_{\mu\raro\mu-1}\bar{\ot}\xi_{\mu-1\raro\mu}\bar{\ot}\xi_{\mu\raro\mu-1},\nonumber \\
    && \xi_{\mu\raro\mu+1}\bar{\ot}\xi_{\mu+1\raro\mu+2}\bar{\ot}\xi_{\mu+2\raro\mu+1},\\
    &&\xi_{\mu\raro\mu-1}\bar{\ot}\xi_{\mu-1\raro\mu-2}\bar{\ot}\xi_{\mu-2\raro\mu-1},\nonumber \\
    && \xi_{\mu\raro\mu+1}\bar{\ot}\xi_{\mu+1\raro\mu}\bar{\ot}\xi_{\mu\raro\mu-1},\nonumber \\
    &&\xi_{\mu\raro\mu-1}\bar{\ot}\xi_{\mu-1\raro\mu}\bar{\ot}\xi_{\mu\raro\mu+1}.\nonumber 
    \end{eqnarray} 
    Note that $\xi_{\mu\raro\mu+1}\bar{\ot}\xi_{\mu+1\raro\mu+2}\bar{\ot}\xi_{\mu+2\raro\mu+3}$ and $\xi_{\mu\raro\mu-1}\bar{\ot}\xi_{\mu-1\raro\mu-2}\bar{\ot}\xi_{\mu-2\raro\mu-3}$ remain invariant under both $\sigma_{12}$ and $\sigma_{23}$. Therefore the braid relation is satisfied trivially on such basis elements. As for other basis elements emitting from a vertex $\mu$, we shall show the braid relation is satisfied for the three length paths between adjacent vertices $\mu$ and $\mu+1$. The braid relations for the paths between the adjacent vertices $\mu$ and $\mu-1$ can be proved similarly.
\begin{align*}
  \sigma_{12}(\xi_{\mu \raro\mu+1}\bar{\ot}\xi_{\mu+1\raro\mu}\bar{\ot}\xi_{\mu\raro\mu+1})&=\xi_{\mu \raro\mu-1}\bar{\ot}\xi_{\mu-1\raro\mu}\bar{\ot}\xi_{\mu\raro\mu+1}\\
    \sigma_{23}\sigma_{12}(\xi_{\mu \raro\mu+1}\bar{\ot}\xi_{\mu+1\raro\mu}\bar{\ot}\xi_{\mu\raro\mu+1})&=\xi_{\mu \raro\mu-1}\bar{\ot}\xi_{\mu-1\raro\mu}\bar{\ot}\xi_{\mu\raro\mu+1}\\
    \sigma_{12}\sigma_{23}\sigma_{12}(\xi_{\mu \raro\mu+1}\bar{\ot}\xi_{\mu+1\raro\mu}\bar{\ot}\xi_{\mu\raro\mu+1})&=\xi_{\mu \raro\mu+1}\bar{\ot}\xi_{\mu+1\raro\mu}\bar{\ot}\xi_{\mu\raro\mu+1}
\end{align*}
On the other hand 
\begin{align*}
    \sigma_{23}(\xi_{\mu \raro\mu+1}\bar{\ot}\xi_{\mu+1\raro\mu}\bar{\ot}\xi_{\mu\raro\mu+1})&=\xi_{\mu \raro\mu+1}\bar{\ot}\xi_{\mu+1\raro\mu+2}\bar{\ot}\xi_{\mu+2\raro\mu+1}\\
     \sigma_{12}\sigma_{23}(\xi_{\mu \raro\mu+1}\bar{\ot}\xi_{\mu+1\raro\mu}\bar{\ot}\xi_{\mu\raro\mu+1})&=\xi_{\mu \raro\mu+1}\bar{\ot}\xi_{\mu+1\raro\mu+2}\bar{\ot}\xi_{\mu+2\raro\mu+1}\\
     \sigma_{23}\sigma_{12}\sigma_{23}(\xi_{\mu \raro\mu+1}\bar{\ot}\xi_{\mu+1\raro\mu}\bar{\ot}\xi_{\mu\raro\mu+1})&=\xi_{\mu \raro\mu+1}\bar{\ot}\xi_{\mu+1\raro\mu}\bar{\ot}\xi_{\mu\raro\mu+1}
\end{align*}
Therefore, 
\begin{align*}
&\sigma_{23}\sigma_{12}\sigma_{23}(\xi_{\mu \raro\mu+1}\bar{\ot}\xi_{\mu+1\raro\mu}\bar{\ot}\xi_{\mu\raro\mu+1})\\
&=\sigma_{12}\sigma_{23}\sigma_{12}(\xi_{\mu \raro\mu+1}\bar{\ot}\xi_{\mu+1\raro\mu}\bar{\ot}\xi_{\mu\raro\mu+1}).
\end{align*}
\indent Similarly 
\begin{align*}
 \sigma_{12}(\xi_{\mu \raro\mu+1}\bar{\ot}\xi_{\mu+1\raro\mu+2}\bar{\ot}\xi_{\mu+2\raro\mu+1})&=\xi_{\mu \raro\mu+1}\bar{\ot}\xi_{\mu+1\raro\mu+2}\bar{\ot}\xi_{\mu+2\raro\mu+1}\\
   \sigma_{23}\sigma_{12}(\xi_{\mu \raro\mu+1}\bar{\ot}\xi_{\mu+1\raro\mu+2}\bar{\ot}\xi_{\mu+2\raro\mu+1})&=\xi_{\mu \raro\mu+1}\bar{\ot}\xi_{\mu+1\raro\mu}\bar{\ot}\xi_{\mu\raro\mu+1}\\
    \sigma_{12}\sigma_{23}\sigma_{12}(\xi_{\mu \raro\mu+1}\bar{\ot}\xi_{\mu+1\raro\mu+2}\bar{\ot}\xi_{\mu+2\raro\mu+1})&=\xi_{\mu \raro\mu-1}\bar{\ot}\xi_{\mu-1\raro\mu}\bar{\ot}\xi_{\mu\raro\mu+1}
   \end{align*}
   \begin{align*}
    \sigma_{23}(\xi_{\mu \raro\mu+1}\bar{\ot}\xi_{\mu+1\raro\mu+2}\bar{\ot}\xi_{\mu+2\raro\mu+1})&=\xi_{\mu \raro\mu+1}\bar{\ot}\xi_{\mu+1\raro\mu}\bar{\ot}\xi_{\mu\raro\mu+1}\\
    \sigma_{12}\sigma_{23}(\xi_{\mu \raro\mu+1}\bar{\ot}\xi_{\mu+1\raro\mu+2}\bar{\ot}\xi_{\mu+2\raro\mu+1})&=\xi_{\mu \raro\mu-1}\bar{\ot}\xi_{\mu-1\raro\mu}\bar{\ot}\xi_{\mu\raro\mu+1}\\
    \sigma_{23}\sigma_{12}\sigma_{23}(\xi_{\mu \raro\mu+1}\bar{\ot}\xi_{\mu+1\raro\mu+2}\bar{\ot}\xi_{\mu+2\raro\mu+1})&=\xi_{\mu \raro\mu-1}\bar{\ot}\xi_{\mu-1\raro\mu}\bar{\ot}\xi_{\mu\raro\mu+1}.
   \end{align*}
   Therefore
   \begin{align*}
   &\sigma_{12}\sigma_{23}\sigma_{12}(\xi_{\mu \raro\mu+1}\bar{\ot}\xi_{\mu+1\raro\mu+2}\bar{\ot}\xi_{\mu+2\raro\mu+1})\\
   &=\sigma_{23}\sigma_{12}\sigma_{23}(\xi_{\mu \raro\mu+1}\bar{\ot}\xi_{\mu+1\raro\mu+2}\bar{\ot}\xi_{\mu+2\raro\mu+1}),
   \end{align*}  And 
   \begin{align*}
    \sigma_{12}(\xi_{\mu \raro\mu-1}\bar{\ot}\xi_{\mu-1\raro\mu}\bar{\ot}\xi_{\mu\raro\mu+1})&=\xi_{\mu \raro\mu+1}\bar{\ot}\xi_{\mu+1\raro\mu}\bar{\ot}\xi_{\mu\raro\mu+1}\\ 
     \sigma_{23}\sigma_{12}(\xi_{\mu \raro\mu-1}\bar{\ot}\xi_{\mu-1\raro\mu}\bar{\ot}\xi_{\mu\raro\mu+1})&=\xi_{\mu \raro\mu+1}\bar{\ot}\xi_{\mu+1\raro\mu+2}\bar{\ot}\xi_{\mu+2\raro\mu+1}\\
    \sigma_{12}\sigma_{23}\sigma_{12}(\xi_{\mu \raro\mu-1}\bar{\ot}\xi_{\mu-1\raro\mu}\bar{\ot}\xi_{\mu\raro\mu+1})&=\xi_{\mu \raro\mu+1}\bar{\ot}\xi_{\mu+1\raro\mu+2}\bar{\ot}\xi_{\mu+2\raro\mu+1}
   \end{align*}
   \begin{align*}
    \sigma_{23}(\xi_{\mu \raro\mu-1}\bar{\ot}\xi_{\mu-1\raro\mu}\bar{\ot}\xi_{\mu\raro\mu+1})&=\xi_{\mu \raro\mu-1}\bar{\ot}\xi_{\mu-1\raro\mu}\bar{\ot}\xi_{\mu\raro\mu+1}\\ 
     \sigma_{12}\sigma_{23}(\xi_{\mu \raro\mu-1}\bar{\ot}\xi_{\mu-1\raro\mu}\bar{\ot}\xi_{\mu\raro\mu+1})&=\xi_{\mu \raro\mu+1}\bar{\ot}\xi_{\mu+1\raro\mu}\bar{\ot}\xi_{\mu\raro\mu+1}\\
    \sigma_{23}\sigma_{12}\sigma_{23}(\xi_{\mu \raro\mu-1}\bar{\ot}\xi_{\mu-1\raro\mu}\bar{\ot}\xi_{\mu\raro\mu+1})&=\xi_{\mu \raro\mu+1}\bar{\ot}\xi_{\mu+1\raro\mu+2}\bar{\ot}\xi_{\mu+2\raro\mu+1}   
   \end{align*}
   So 
   \begin{align*}
     &\sigma_{12}\sigma_{23}\sigma_{12}(\xi_{\mu \raro\mu-1}\bar{\ot}\xi_{\mu-1\raro\mu}\bar{\ot}\xi_{\mu\raro\mu+1})\\
     &=\sigma_{23}\sigma_{12}\sigma_{23}(\xi_{\mu \raro\mu-1}\bar{\ot}\xi_{\mu-1\raro\mu}\bar{\ot}\xi_{\mu\raro\mu+1})  
   \end{align*}
   Therefore the braid relation holds on the $\mathbb{C}$-linear basis elements of $\Omega^{1}(\cla)^{\bar{\otimes}3}$. By obvious $\mathbb{C}$-linearity of the maps $\sigma_{12},\sigma_{23}$, the braid relation holds on $\Omega^{1}(\cla)^{\bar{\otimes}3}$.
\end{proof}
Thanks to the above Lemma \ref{sigmaok}, by Theorem \ref{prolongationgraph}, we can construct the canonical prolongation of the first order $\ast$-calculus with respect to $\sigma$ to obtain a $\ast$-differential calculus.
\begin{theorem}
    \label{polygoncalculus}
    The canonical prolongation of the first order $\ast$-differential calculus on the bidirected $n$-gon with respect to the bimodule map $\sigma$ is an orientable $\ast$-differential calculus of dimension $2$. Moreover, $\Omega^{1}(\cla)$ and $\Omega^{2}(\cla)$ are both left and right parallelizable in the sense of \cite{Majid}.
\end{theorem}
\begin{proof}
    For the bimodule of $2$-forms, observe that the antisymmetrization map $A_{2}$ is $(\sigma-{\rm Id})$ on $\Omega^{1}(\cla)\bar{\ot}\Omega^{1}(\cla)$. The first two types of the basis elements of $\Omega^{1}(\cla)\bar{\ot}\Omega^{1}(\cla)$ as in formula (\ref{omega2basis}) are fixed by $\sigma$ and therefore mapped to zero elements in the quotient bimodule $\Omega^{2}(\cla)$. As for the rest of the basis elements, it follows from the definition of $\sigma$ that $(\xi_{\mu\raro\mu+1}\bar{\ot}\xi_{\mu+1\raro\mu}+\xi_{\mu\raro\mu-1}\bar{\ot}\xi_{\mu-1\raro\mu})$ belong to the kernel of $A_{2}$ for all $\mu=1,\ldots, n$ and therefore $\xi_{\mu\raro\mu+1}\wedge\xi_{\mu+1\raro\mu}=-\xi_{\mu\raro\mu-1}\wedge\xi_{\mu-1\raro\mu}$ in $\Omega^{2}(\cla)$. So $\Omega^{2}(\cla)$ has a $\mathbb{C}$-linear basis $\{\xi_{\mu\raro\mu-1}\wedge\xi_{\mu-1\raro\mu}\}_{\mu=1,\ldots, n}$ or $\{\xi_{\mu\raro\mu+1}\wedge\xi_{\mu+1\raro\mu}\}_{\mu=1,\ldots, n}$.
    
    Now to calculate $\Omega^{3}(\cla)$, let us write down the antisymmetrization map $A_{3}$ on $\Omega^{1}(\cla)^{\bar{\otimes}3}$ explicitly. To do this, we write all the elements in $S_{3}$ as a product of neighboring transpositions:
    \begin{eqnarray*}
        (1 \ 2 \ 3)=(1\ 2)(2\ 3), \ (1 \ 3 \ 2)=(2\ 3)(1\ 2), \ (1\ 3)=(1\ 2)(2\ 3)(1\ 2).
    \end{eqnarray*}Therefore, by definition,
    \begin{displaymath}
        A_{3}\equiv{\rm Id}+\sigma_{12}\sigma_{23}+\sigma_{23}\sigma_{12}-\sigma_{12}-\sigma_{23}-\sigma_{12}\sigma_{23}\sigma_{12}.
    \end{displaymath}
    Recall the collection of basis elements of $\Omega^{1}(\cla)^{\bar{\otimes}3}$ emitting from a vertex $\mu$ as given in the proof of the Lemma \ref{sigmaok}. As observed there,
    the first two three-length paths in formula (\ref{omega3basis}) are invariant under both $\sigma_{12},\sigma_{23}$ and therefore belong to the kernel of $A_{3}$. As for the rest of the paths, we prove that the three length paths between the adjacent vertices $\mu,\mu+1$ belong to the kernel. The rest of the three-length paths between adjacent vertices $\mu,\mu-1$ can be shown to belong to the kernel of $A_{3}$ similarly. Borrowing the calculations from the proof of the Lemma \ref{sigmaok}, we have
    \begin{eqnarray*}
      \begin{aligned}
        &A_{3}(\xi_{\mu\raro\mu+1}\bar{\ot}\xi_{\mu+1\raro\mu}\bar{\ot}\xi_{\mu\raro\mu+1})\\
        &= ({\rm Id}-\sigma_{12}-\sigma_{23}-\sigma_{12}\sigma_{23}\sigma_{12}+\sigma_{12}\sigma_{23}+\sigma_{23}\sigma_{12})(\xi_{\mu\raro\mu+1}\bar{\ot}\xi_{\mu+1\raro\mu}\bar{\ot}\xi_{\mu\raro\mu+1})\\
        &= \xi_{\mu\raro\mu+1}\bar{\ot}\xi_{\mu+1\raro\mu}\bar{\ot}\xi_{\mu\raro\mu+1}-\xi_{\mu\raro\mu-1}\bar{\ot}\xi_{\mu-1\raro\mu}\bar{\ot}\xi_{\mu\raro\mu+1}\\
        &\quad-\xi_{\mu\raro\mu+1}\bar{\ot}\xi_{\mu+1\raro\mu+2}\bar{\ot}\xi_{\mu+2\raro\mu+1}-\xi_{\mu\raro\mu+1}\bar{\ot}\xi_{\mu+1\raro\mu}\bar{\ot}\xi_{\mu\raro\mu+1}\\
        &\quad+ \xi_{\mu\raro\mu+1}\bar{\ot}\xi_{\mu+1\raro\mu+2}\bar{\ot}\xi_{\mu+2\raro\mu+1}+\xi_{\mu\raro\mu-1}\bar{\ot}\xi_{\mu-1\raro\mu}\bar{\ot}\xi_{\mu\raro\mu+1}\\
        &= 0.
      \end{aligned}
    \end{eqnarray*}
    Similarly,
    \begin{eqnarray*}
      \begin{aligned}
       & A_{3}(\xi_{\mu\raro\mu+1}\bar{\ot}\xi_{\mu+1\raro\mu+2}\bar{\ot}\xi_{\mu+2\raro\mu+1})\\
        &= ({\rm Id}-\sigma_{12}-\sigma_{23}-\sigma_{12}\sigma_{23}\sigma_{12}+\sigma_{12}\sigma_{23}+\sigma_{23}\sigma_{12})(\xi_{\mu\raro\mu+1}\bar{\ot}\xi_{\mu+1\raro\mu+2}\bar{\ot}\xi_{\mu+2\raro\mu+1})\\
        &=\xi_{\mu\raro\mu+1}\bar{\ot}\xi_{\mu+1\raro\mu+2}\bar{\ot}\xi_{\mu+2\raro\mu+1}-\xi_{\mu\raro\mu+1}\bar{\ot}\xi_{\mu+1\raro\mu+2}\bar{\ot}\xi_{\mu+2\raro\mu+1}\\
        &\quad- \xi_{\mu\raro\mu+1}\bar{\ot}\xi_{\mu+1\raro\mu}\bar{\ot}\xi_{\mu\raro\mu+1}-\xi_{\mu\raro\mu-1}\bar{\ot}\xi_{\mu-1\raro\mu}\bar{\ot}\xi_{\mu\raro\mu+1}\\
        &\quad+ \xi_{\mu\raro\mu-1}\bar{\ot}\xi_{\mu-1\raro\mu}\bar{\ot}\xi_{\mu\raro\mu+1}+ \xi_{\mu\raro\mu+1}\bar{\ot}\xi_{\mu+1\raro\mu}\bar{\ot}\xi_{\mu\raro\mu+1}\\
        &= 0,
      \end{aligned}
    \end{eqnarray*}
    and
    \begin{eqnarray*}
      \begin{aligned}
      &A_{3}(\xi_{\mu\raro\mu-1}\bar{\ot}\xi_{\mu-1\raro\mu}\bar{\ot}\xi_{\mu\raro\mu+1})\\
      &= ({\rm Id}-\sigma_{12}-\sigma_{23}-\sigma_{12}\sigma_{23}\sigma_{12}+\sigma_{12}\sigma_{23}+\sigma_{23}\sigma_{12})(\xi_{\mu\raro\mu-1}\bar{\ot}\xi_{\mu-1\raro\mu}\bar{\ot}\xi_{\mu\raro\mu+1})\\
      &= \xi_{\mu\raro\mu-1}\bar{\ot}\xi_{\mu-1\raro\mu}\bar{\ot}\xi_{\mu\raro\mu+1}-\xi_{\mu\raro\mu+1}\bar{\ot}\xi_{\mu+1\raro\mu}\bar{\ot}\xi_{\mu\raro\mu+1}\\
      &\quad- \xi_{\mu\raro\mu-1}\bar{\ot}\xi_{\mu-1\raro\mu}\bar{\ot}\xi_{\mu\raro\mu+1}- \xi_{\mu\raro\mu+1}\bar{\ot}\xi_{\mu+1\raro\mu+2}\bar{\ot}\xi_{\mu+2\raro\mu+1}\\ &\quad+\xi_{\mu\raro\mu+1}\bar{\ot}\xi_{\mu+1\raro\mu}\bar{\ot}\xi_{\mu\raro\mu+1}+\xi_{\mu\raro\mu+1}\bar{\ot}\xi_{\mu+1\raro\mu+2}\bar{\ot}\xi_{\mu+2\raro\mu+1}\\
      &= 0.
      \end{aligned}
    \end{eqnarray*}
        Hence all the basis elements of $\Omega^{1}(\cla)^{\bar{\otimes}3}$ are in the kernel of the antisymmetrization map $A_{3}$ proving that $\Omega^{3}(\cla)=0$. By the remark before Lemma \ref{star_2form} , $\Omega^{n}(\cla)=0$ for all $n\geq 3$.
    This proves that the dimension of the calculus is $2$.
    
     $\Omega^{1}(\cla)$ is a free module of rank $2$ (both left and right module over $\cla$). This follows from the fact that each vertex of the graph $n$-gon receives and emits $2$ edges. In fact, following the argument as in \cite{Majid}, $\{\sum_{\mu=1}^{n}\xi_{\mu\raro\mu+1},\sum_{\mu=1}^{n}\xi_{\mu\raro\mu-1}\}$ is an $\cla$-linear basis both as right and left module. As for $\Omega^{2}(\cla)$, it is easy to see that the right and left actions of elements of $\cla$ on $\Omega^{2}(\cla)$ coincide and $\{\sum_{\mu=1}^{n}\xi_{\mu\raro\mu+1}\wedge\xi_{\mu+1\raro\mu}\}$ is an $\cla$-linear basis both as right and left module. This proves that $\Omega^{i}(\cla)$ are both left and right parallelizable in the sense of \cite{Majid}.
     
     Consider the map $\Phi:\cla\raro\Omega^{2}(\cla)$ given by
    \begin{displaymath}        
    \Phi(f)=\sum_{\mu=1}^{n}f(\mu)\xi_{\mu\raro\mu-1}\wedge\xi_{\mu-1\raro\mu}.
    \end{displaymath}
Using the fact that the right and left action of $ 
\cla$ on $\Omega^{2}(\cla)$ coincide, clearly $\Phi$ is an $\cla$\ndash$\cla$-bimodule isomorphism proving that the calculus is orientable.  
\end{proof}
\begin{remark}
    The isomorphism $\Phi$ corresponds to the choice of $\cla$-basis $\{\sum_{\mu=1}^{n}\xi_{\mu\raro\mu-1}\wedge\xi_{\mu-1\raro\mu}\}$ of $\Omega^{2}(\cla)$. One could choose the bimodule isomorphism $\Phi$ with respect to the choice of $\cla$-basis $\{\sum_{\mu=1}^{n}\xi_{\mu\raro\mu+1}\wedge\xi_{\mu+1\raro\mu}\}$ of $\Omega^{2}(\cla)$. That amounts to choosing an `opposite orientation'. Our choice of orientation is driven by positivity of a certain Hochschild cocycle which will be discussed in the last section.
\end{remark}
We can define a state $\tau$ on $C(V)$ by $\tau(f):=\sum_{\mu=1}^{n}\frac{1}{n}f(\mu)$. Then with this state at our disposal we have the following:
\begin{lemma}
    \label{gradedtrace}
    Let $(\Omega^{\bullet}(\cla),\der,\ast)$ be the $2$-dimensional calculus in Theorem \ref{polygoncalculus}. Then there is a closed graded trace of dimension $2$ on $(\Omega^{\bullet}(\cla),\der,\ast)$ given by 
    \begin{displaymath}
        \int \omega:=\tau(\Phi^{-1}(\omega)),
    \end{displaymath}
    where $\Phi:\cla\raro\Omega^{2}(\cla)$ is the bimodule isomorphism in Theorem \ref{polygoncalculus}.
\end{lemma}
\begin{proof}
    To prove the equality $\int \der\omega=0$ for all $\omega\in\Omega^{1}(\cla)$, by $\mathbb{C}$-linearity of $\int$ and $\der$, it is enough to show that $\int \der\omega=0$ where $\omega$ is a $\mathbb{C}$-linear basis element of $\Omega^{1}(\cla)$. For $\mu=1,\ldots,n$, 
    \begin{eqnarray*}
      \begin{aligned}
        \der(\xi_{\mu\raro\mu+1})&=\xi_{\mu\raro\mu+1}\wedge\theta+\theta\wedge\xi_{\mu\raro\mu+1}\\
        &=\xi_{\mu\raro\mu+1}\wedge\xi_{\mu+1\raro\mu}+\xi_{\mu+1\raro\mu}\wedge\xi_{\mu\raro\mu+1}\\
        &=-\xi_{\mu\raro\mu-1}\wedge\xi_{\mu-1\raro\mu}+\xi_{\mu+1\raro\mu}\wedge\xi_{\mu\raro\mu+1}
      \end{aligned}
    \end{eqnarray*}
    Therefore, for $\mu=1,\ldots, n$, $\Phi^{-1}(\der\xi_{\mu\raro\mu+1})=f\in\cla$ where $f(\mu)=-1,f(\mu+1)=1$ and $f$ is zero on the rest of the vertices. Hence by definition, $\int \der\xi_{\mu\raro\mu+1}=\frac{1}{n}(1-1)=0$. Similarly one can show that for the rest of the basis elements also $\int \der\omega=0$. \\
    \indent For the second condition of Definition \ref{defgradedtrace}, take $\omega_{1}, \omega_{2}\in\Omega^{1}(\cla)$. Then there are scalars $\{\alpha_{\mu},\alpha^{\mu},\beta_{\mu},\beta^{\mu}\}_{\mu=1,\ldots,n}$ such that
    \begin{eqnarray*}
        &&\omega_{1}=\sum_{\mu=1}^{n}\bigl(\alpha_{\mu}\xi_{\mu\raro\mu+1}+\alpha^{\mu}\xi_{\mu\raro\mu-1}\bigr),\\
        &&\omega_{2}=\sum_{\mu=1}^{n}\bigl(\beta_{\mu}\xi_{\mu\raro\mu+1}+\beta^{\mu}\xi_{\mu\raro\mu-1}\bigr).
    \end{eqnarray*}
    Therefore 
\begin{eqnarray*}
  \begin{aligned}
\omega_{1}\wedge\omega_{2}+\omega_{2}\wedge\omega_{1}
&=\sum_{\mu=1}^{n}(\alpha^{\mu}\beta_{\mu-1}-\alpha_{\mu}\beta^{\mu+1})\xi_{\mu\raro\mu-1}\wedge\xi_{\mu-1\raro\mu}\\
&\qquad+\sum_{\mu=1}^{n}(\beta^{\mu}\alpha_{\mu-1}-\beta_{\mu}\alpha^{\mu+1})\xi_{\mu\raro\mu-1}\wedge\xi_{\mu-1\raro\mu}.
\end{aligned}
    \end{eqnarray*}
    By the definition of $\int$,
    \begin{equation*}        \int(\omega_{1}\wedge\omega_{2}+\omega_{2}\wedge\omega_{1})
        =\sum_{\mu=1}^{n}\bigl(\alpha^{\mu}\beta_{\mu-1}-\alpha_{\mu}\beta^{\mu+1}+\beta^{\mu}\alpha_{\mu-1}-\beta_{\mu}\alpha^{\mu+1}\bigr)
        =0
    \end{equation*}
    This proves the second condition of Definition \ref{defgradedtrace}.
\end{proof}
Now we shall discuss an almost complex structure on $(\Omega^{\bullet}(\cla),\der,\ast)$. We start by defining a bimodule map $J:\Omega^{1}(\cla)\raro\Omega^{1}(\cla)$ satisfying $J^2=-{\rm Id}$ and $J(\omega^{\ast})=(J(\omega))^{\ast}$. As usual we define $J$ on the basis elements of $\Omega^{1}(\cla)$ and extend it $\mathbb{C}$- linearly. 
\begin{eqnarray*}
    J(\xi_{\mu\raro\mu+1})=i\xi_{\mu\raro\mu+1}, \ \mu=1,2,\ldots,n
\end{eqnarray*}
As we have already observed compatibility of $J$ with $\ast$ fixes the $-i$-th eigenspace of $J$ i.e. $J$ is forced to be defined on the rest of the basis elements by the following:
\begin{eqnarray*}
J(\xi_{\mu\raro\mu-1})=-i\xi_{\mu\raro\mu-1},\ \mu=1,2,\ldots,n.
\end{eqnarray*}
With the above definition, it is trivial to check that $J$ is a bimodule map; $J^{2}=-{\rm Id}$ and $(J(\omega))^{\ast}=J(\omega^{\ast})$. Moreover, $\Omega^{1}(\cla)=\Omega^{1,0}(\cla)\oplus\Omega^{0,1}(\cla)$, where
\begin{eqnarray*}
&&\Omega^{1,0}(\cla)={\rm Sp}_{\mathbb{C}}\{\xi_{1\raro 2},\xi_{2\raro 3},\ldots,\xi_{n\raro 1}\},\\
&&\Omega^{0,1}(\cla)={\rm Sp}_{\mathbb{C}}\{\xi_{1\raro n},\xi_{n\raro n-1},\ldots,\xi_{2\raro 1}\}.
\end{eqnarray*}
Now by Theorem \ref{extension}, $J$ will extend to $(\Omega^{\bullet}(\cla),\der,\ast)$ as a degree zero derivation to define an integrable complex structure if the extension of $J$ on $\Omega^{1}(\cla)\bar{\ot}\Omega^{1}(\cla)$ given by the formula \eqref{Jextension}, commutes with $\sigma$. 
\begin{lemma}
\label{lemmaJcommutsigma}
    The extension of $J$ on $\Omega^{1}(\cla)\bar{\ot}\Omega^{1}(\cla)$ commutes with the bimodule map $\sigma$.
\end{lemma}
\begin{proof}
    Again as before it is enough to prove commutation of $J$ and $\sigma$ on the basis elements of $\Omega^{1}(\cla)\bar{\ot}\Omega^{1}(\cla)$. The first two types of basis elements of $\Omega^{1}(\cla)\bar{\ot}\Omega^{1}(\cla)$ in the formula \eqref{omega2basis} are in the eigen spaces of $J$ corresponding to the eigen values $\pm 2i$ and are fixed by $\sigma$. The rest of the basis elements are in the eigen space corresponding to the eigenvalue $0$. The basis elements in the $0$ eigenspace are again mapped to the basis elements in the $0$ eigenspace. This can easily be seen from the formula of $\sigma$. Therefore $J$ commutes with $\sigma$.
\end{proof}
We write down the easy to see formula for the bimodule maps $\pi^{0,1}$ and $\pi^{1,0}$ on the $\mathbb{C}$-linear basis elements of $\Omega^{1}(\cla)$ which will be used later:
\begin{equation}
\label{piformula}
  \begin{aligned}
    &\pi^{0,1}(\xi_{\mu\raro\mu+1})=0;  &\pi^{0,1}(\xi_{\mu\raro\mu-1})=\xi_{\mu\raro\mu-1};\\
    &\pi^{1,0}(\xi_{\mu\raro\mu+1})=\xi_{\mu\raro\mu+1}; &\pi^{1,0}(\xi_{\mu\raro\mu-1})=0.
  \end{aligned}
\end{equation}
As $\Omega^{k}(\cla)=0$ for $k\geq 3$, there are only $(0,1)$, $(1,0)$, $(2,0)$, $(0,2)$ and $(1,1)$ forms. From the proof of Lemma \ref{lemmaJcommutsigma}, it is easy to see that $\Omega^{2,0}(\cla)=0=\Omega^{0,2}$ and therefore $\Omega^{2}(\cla)=\Omega^{1,1}(\cla)$. We can consider the following Dolbeault complex
\begin{align}
\label{Dolbeault}
\begin{array}{ccccccc}
  0 &\longrightarrow & \Omega^{1,0}(\cla) & \overset{\overline{\partial}}{\longrightarrow} & \Omega^{1,1}(\cla) & \overset{\overline{\partial}}{\longrightarrow} & 0 
  \end{array}
  \end{align}
to compute the Dolbeault cohomomology vector spaces $H^{1,0}_{\overline{\partial}}(\cla)$ and $H^{1,1}_{\overline{\partial}}(\cla)$. Rest of the Dolbeault cohomology spaces vanish trivially. 
\begin{lemma}
    The Dolbeault cohomology of the Dolbeault complex \eqref{Dolbeault} is given by
    \[
        H^{1,0}_{\overline{\partial}}(\cla)=\mathbb{C} =  H^{1,1}_{\overline{\partial}}(\cla),
        \quad
        H^{1,q}_{\overline{\partial}}(\cla)=0 \text{ for all } q\geq 2.
    \]
\end{lemma}
\begin{proof}
It is clear that $H^{1,0}_{\overline{\partial}}(\cla)={\rm Ker}(\overline{\partial})$ and $H^{1,1}_{\overline{\partial}}(\cla)=\Omega^{1,1}(\cla)/{\rm Im}(\overline{\partial})$. A $\mathbb{C}$-linear basis of $\Omega^{1,0}(\cla)$ is given by $\{\xi_{\mu\raro\mu+1}\}_{\mu=1,\ldots,n}$. By definition,
\begin{eqnarray*}
  \begin{aligned}
    \overline{\partial}(\xi_{\mu\raro\mu+1})
    &=\pi^{1,1}(\der\xi_{\mu\raro\mu+1})\\
    &= \pi^{1,1}(\theta\wedge\xi_{\mu\raro\mu+1}+\xi_{\mu\raro\mu+1}\wedge\theta)\\
    &= \xi_{\mu\raro\mu+1}\wedge\xi_{\mu+1\raro\mu}+\xi_{\mu+1\raro\mu}\wedge\xi_{\mu\raro\mu+1}\\
    &= \xi_{\mu\raro\mu+1}\wedge\xi_{\mu+1\raro\mu}-\xi_{\mu+1\raro\mu+2}\wedge\xi_{\mu+2\raro\mu+1}.
  \end{aligned}
\end{eqnarray*}
Therefore an element $\sum_{\mu=1}^{n}\alpha_{\mu}\xi_{\mu\raro\mu+1}\in{\rm Ker}(\overline{\partial}) \ (\alpha_{\mu}\in\mathbb{C})$ implies that
\begin{eqnarray*}
    &&\sum_{\mu=1}^{n}\bigl(\alpha_{\mu}\xi_{\mu\raro\mu+1}\wedge\xi_{\mu+1\raro\mu}-\alpha_{\mu}\xi_{\mu+1\raro\mu+2}\wedge\xi_{\mu+2\raro\mu+1}\bigr)=0\\
    &\Rightarrow&\sum_{\mu=1}^{n}(\alpha_{\mu}-\alpha_{\mu-1})\xi_{\mu\raro\mu+1}\wedge\xi_{\mu+1\raro\mu}=0
\end{eqnarray*}
Hence, $\alpha_{\mu}=\alpha_{\mu-1}$ for all $\mu=1,\ldots,n$ i.e $\alpha_{1}=\alpha_{2}=\cdots=\alpha_{n}$. So 
\begin{displaymath}
    H^{1,0}_{\overline{\partial}}={\rm Sp}_{\mathbb{C}}\{\sum_{\mu=1}^{n}\xi_{\mu\raro\mu+1}\}\cong\mathbb{C}.
\end{displaymath}
To determine $H^{1,1}_{\overline{\partial}}(\cla)$, note that by rank nullity theorem applied to the linear map $\overline{\partial}:\Omega^{1,0}(\cla)\raro\Omega^{1,1}(\cla)$, ${\rm dim}_{\mathbb{C}}({\rm Im}\overline{\partial})=n-1$. Therefore ${\rm dim}_{\mathbb{C}}\Omega^{1,1}(\cla)$ being $n$, ${\rm dim}_{\mathbb{C}}(H^{1,1}_{\overline{\partial}}(\cla))=1$ and therefore $H^{1,1}_{\overline{\partial}}(\cla)\cong\mathbb{C}$. Class of any one of the basis elements $\{\xi_{\mu\raro\mu+1}\wedge\xi_{\mu+1\raro\mu}\}_{\mu=1}^{n}$ in $H^{1,1}_{\overline{\partial}}(\cla)$ gives a basis. The elements $\{\xi_{\mu\raro\mu+1}\wedge\xi_{\mu+1\raro\mu}\}_{\mu=1}^{n}$ are cohomologous to each other which is clear from the expression of $\{\overline{\partial}(\xi_{\mu\raro\mu+1})\}_{\mu=1}^{n}$.
\end{proof}
\subsection{A Holomorphic structure on the exterior bundle} In this section we will explore a holomorphic structure on the `exterior bundle'. This `bundle' is obtained through the canonical prolongation of the first order diffrential calculus originating from a bidirected polygon structure with respect to the bimodule map $\sigma$ of permutation type. First we shall induce a holomorphic structure on $\Omega^{1}(\cla)$ coming from a bimodule connection $\nabla$ on $\Omega^{1}(\cla)$ with respect to the same bimodule isomorphism $\sigma:\Omega^{1}(\cla)\bar{\ot}\Omega^{1}(\cla)\raro\Omega^{1}(\cla)\bar{\ot}\Omega^{1}(\cla)$ used in the prolongation procedure. Note that since $\Omega^{0,2}(\cla)=0$, any such connection inherently induces a holomorphic structure on $\Omega^{1}(\cla)$. 
\begin{lemma}
    \label{bimodisomorphrestrcit}
    The map $\sigma$ satisfies $(\pi^{1,0}\bar{\ot}{\rm id})\circ\sigma\circ({\rm id}\bar{\ot}\pi^{0,1})=0=(\pi^{0,1}\bar{\ot}{\rm id})\circ\sigma\circ({\rm id}\bar{\ot}\pi^{1,0})$ on $\Omega^{1}(\cla)\bar{\ot}\Omega^{1}(\cla)$. In particular $\sigma$ restricts to a bimodule isomorphism between $\Omega^{1}(\cla)\bar{\ot}\Omega^{0,1}(\cla)$ and $\Omega^{0,1}(\cla)\bar{\ot}\Omega^{1}(\cla)$; $\sigma\circ({\rm id}\bar{\ot}\pi^{0,1})=(\pi^{0,1}\bar{\ot}{\rm id})\circ\sigma$ on $\Omega^{1}(\cla)\bar{\ot}\Omega^{1}(\cla)$. 
\end{lemma}
\begin{proof}
  Again it is enough to show the identities on the $\mathbb{C}$-linear basis elements of $\Omega^{1}(\cla)\bar{\ot}\Omega^{1}(\cla)$. Recalling the basis elements from the formula (\ref{omega2basis}), we do the following calculations using the formulas \eqref{piformula}.   
  \begin{eqnarray*}
    \begin{aligned}
      &(\pi^{1,0}\bar{\ot}{\rm id})\circ\sigma\circ({\rm id}\bar{\ot}\pi^{0,1})(\xi_{\mu\raro\mu+1}\bar{\ot}\xi_{\mu+1\raro\mu})\\
      &\quad= (\pi^{1,0}\bar{\ot}\rm{\id})\circ\sigma(\xi_{\mu\raro\mu+1}\bar{\ot}\xi_{\mu+1\raro\mu})\\
      &\quad= (\pi^{1,0}\bar{\ot}\rm{\id})(\xi_{\mu\raro\mu-1}\bar{\ot}\xi_{\mu-1\raro\mu})\\
      &\quad= 0.
    \end{aligned}
  \end{eqnarray*}
  $((\pi^{1,0}\bar{\ot}{\rm id})\circ\sigma\circ({\rm id}\bar{\ot}\pi^{0,1}))(\xi_{\mu\raro\mu+1}\bar{\ot}\xi_{\mu+1\raro\mu+2})$ is trivially zero. As for the other basis elements,
  \begin{eqnarray*}
    \begin{aligned}
      &(\pi^{1,0}\bar{\ot}{\rm id})\circ\sigma\circ({\rm id}\bar{\ot}\pi^{0,1})(\xi_{\mu\raro\mu-1}\bar{\ot}\xi_{\mu-1\raro\mu-2})\\
      &\quad=(\pi^{1,0}\bar{\ot}\rm{\id})\circ\sigma(\xi_{\mu\raro\mu-1}\bar{\ot}\xi_{\mu-1\raro\mu-2})\\
      &\quad= (\pi^{1,0}\bar{\ot}\rm{\id})(\xi_{\mu\raro\mu-1}\bar{\ot}\xi_{\mu-1\raro\mu-2})\\
      &\quad= 0.
    \end{aligned}
  \end{eqnarray*}
  Again $((\pi^{1,0}\bar{\ot}{\rm id})\circ\sigma\circ({\rm id}\bar{\ot}\pi^{0,1}))(\xi_{\mu\raro\mu-1}\bar{\ot}\xi_{\mu-1\raro\mu})$ is trivially zero.
\end{proof}
We denote the restriction of $\sigma$ on $\Omega^{1}(\cla)\bar{\ot}\Omega^{0,1}(\cla)$ by $\overline{\sigma}$. Now we shall consider a bimodule connection $\nabla$ on the bimodule $\Omega^{1}(\cla)$ with respect to the bimodule map $\sigma$ used in the canonical prolongation and $\alpha:\Omega^{1}(\cla)\raro\Omega^{1}(\cla)\bar{\ot}\Omega^{1}(\cla)$ as the zero map as in \cite{Majid}. Recall that the formula of such a bimodule connection on $\Omega^{1}(\cla)$ is given by
\begin{eqnarray}
\label{nablaformula}
    \nabla(\omega):=\theta\bar{\ot}\omega-\sigma(\omega\bar{\ot}\theta).
\end{eqnarray}
Recalling that $\theta=\sum_{\mu}(\xi_{\mu\raro\mu+1}+\xi_{\mu\raro\mu-1})$, it is easy to see from the formula (\ref{nablaformula}) that $\nabla$ is given on the basis elements of $\Omega^{1}(\cla)$ by
\begin{eqnarray*}
    \begin{aligned}  
    \nabla(\xi_{\mu\raro\mu+1})=\xi_{\mu-1\raro\mu}&\bar{\ot}\xi_{\mu\raro\mu+1}+\xi_{\mu+1\raro\mu}\bar{\ot}\xi_{\mu\raro\mu+1}\\
    &-\xi_{\mu\raro\mu-1}\bar{\ot}\xi_{\mu-1\raro\mu}
    -\xi_{\mu\raro\mu+1}\bar{\ot}\xi_{\mu+1\raro\mu+2};
    \end{aligned}
    \end{eqnarray*}
    \begin{eqnarray*}
    \begin{aligned}
    \nabla(\xi_{\mu\raro\mu-1})=\xi_{\mu+1\raro\mu}&\bar{\ot}\xi_{\mu\raro\mu-1}+\xi_{\mu-1\raro\mu}\bar{\ot}\xi_{\mu\raro\mu-1}\\
    &-\xi_{\mu\raro\mu-1}\bar{\ot}\xi_{\mu-1\raro\mu-2}
    -\xi_{\mu\raro\mu+1}\bar{\ot}\xi_{\mu+1\raro\mu}.
    \end{aligned}
\end{eqnarray*}
\begin{remark}
    We do not compute the curvature $\nabla^2$ of the connection $\nabla$ here. But it can be shown that $\nabla^2$ is non zero. For example a simple computation reveals that $\nabla^2$ acts on a basis element $\xi_{\mu\raro\mu+1}$ in the following way:
    \begin{displaymath}
        \nabla^{2}(\xi_{\mu\raro\mu+1})=\xi_{\mu-1\raro\mu}\wedge\xi_{\mu\raro\mu-1}\bar{\ot}\xi_{\mu-1\raro\mu}+\xi_{\mu+1\raro\mu}\wedge\xi_{\mu\raro\mu+1}\bar{\ot}\xi_{\mu+1\raro\mu+2}
    \end{displaymath}
    which is clearly non zero.
\end{remark}
As mentioned earlier the above connection will induce a holomorphic structure on $\Omega^{1}(\cla)$ since $\Omega^{0,2}(\cla)=0$. Noting that for all $\mu$, $\pi^{0,1}(\xi_{\mu\raro\mu+1})=0$ and $\pi^{0,1}(\xi_{\mu\raro\mu-1})=\xi_{\mu\raro\mu-1}$, the corresponding $\overline{\partial}$-operator $\overline{\nabla}:=(\pi^{0,1}\bar{\ot}{\rm id})\circ\nabla$ acts on the basis elements of $\Omega^{1}(\cla)$ by the following formulae:
\begin{eqnarray}
    \label{nablabarformula}
    && \overline{\nabla}(\xi_{\mu\raro\mu+1})=\xi_{\mu+1\raro\mu}\bar{\ot}\xi_{\mu\raro\mu+1}-\xi_{\mu\raro\mu-1}\bar{\ot}\xi_{\mu-1\raro\mu};\\
    && \overline{\nabla}(\xi_{\mu\raro\mu-1})=\xi_{\mu+1\raro\mu}\bar{\ot}\xi_{\mu\raro\mu-1}-\xi_{\mu\raro\mu-1}\bar{\ot}\xi_{\mu-1\raro\mu-2}.
\end{eqnarray}
\begin{lemma}
    $\overline{\nabla}$ is a bimodule $\overline{\partial}$ operator on $\Omega^{1}(\cla)$ with respect to the bimodule isomorphism $\overline{\sigma}$ i.e. $\overline{\nabla}$ is a left $\overline{\partial}$-operator satisfying:
    \begin{displaymath}
    \overline{\nabla}(\omega.f)=\overline{\nabla}(\omega).f+\overline{\sigma}(\omega\bar{\ot}\overline{\partial}f).
    \end{displaymath}
\end{lemma}
\begin{proof}
We shall only prove the twisted right Leibniz rule. The left Liebniz rule can be similarly proved and will be left to the reader. 
\begin{eqnarray*}
\begin{aligned}
\overline{\nabla}(\omega.f)&=(\pi^{0,1}\bar{\ot}{\rm id})\circ\nabla(\omega.f)\\
&= (\pi^{0,1}\bar{\ot}{\rm id})(\nabla(\omega).f+\sigma(\omega\bar{\ot}\der{f}))\\
&= \bigl((\pi^{0,1}\bar{\ot}{\rm id})\circ\nabla(\omega)\bigr).f+(\pi^{0,1}\bar{\ot}{\rm id})\circ\sigma(\omega\bar{\ot}\der{f})\\
&= \bigl((\pi^{0,1}\bar{\ot}{\rm id})\circ\nabla(\omega)\bigr).f+\sigma\circ({\rm id}\bar{\ot}\pi^{0,1})(\omega\bar{\ot}\der{f}) \text{ (by Lemma~\ref{bimodisomorphrestrcit})}\\
&= \overline{\nabla}(\omega).f+\overline{\sigma}(\omega\bar{\ot}\overline{\partial}f)).
\end{aligned}
\end{eqnarray*}
\end{proof}

Therefore from the discussion after Definition \ref{bimoduleconnection}, it follows that $\overline{\nabla}^{(\otimes2)}:=(\overline{\nabla}\bar{\ot}1)+\overline{\sigma}_{12}(1\bar{\ot}\overline{\nabla})$ is again a bimodule $\overline{\partial}$-operator with respect to the bimodule isomorphism $\overline{\sigma}_{12}\overline{\sigma}_{23}:\Omega^{1}(\cla)\bar{\ot}\Omega^{1}(\cla)\bar{\ot}\Omega^{0,1}(\cla)\raro\Omega^{0,1}(\cla)\bar{\ot}\Omega^{1}(\cla)\bar{\ot}\Omega^{1}(\cla) $, whereas $\nabla^{(\otimes2)}:=(\nabla\bar{\ot}1)+\sigma_{12}(1\bar{\ot}\nabla)$ is a bimodule connection on $\Omega^{1}(\cla)\bar{\ot}\Omega^{1}(\cla)$ with respect to the bimodule isomorphism $\sigma_{12}\sigma_{23}$. 
\begin{lemma}
    \label{relationnablanablabar}
    \begin{eqnarray}
        \overline{\nabla}^{(\otimes2)}=(\pi^{0,1}\bar{\ot}{\rm id}\bar{\ot}{\rm id})\circ\nabla^{(\otimes2)}.
    \end{eqnarray}
\end{lemma}
\begin{proof}
    \begin{eqnarray*}
      \begin{aligned}
        (\pi^{0,1}\bar{\ot}{\rm id}\bar{\ot}{\rm id})\circ\nabla^{(\otimes2)}&= (\pi^{0,1}\bar{\ot}{\rm id}\bar{\ot}{\rm id})((\nabla\bar{\ot}1)+\sigma_{12}(1\bar{\ot}\nabla))\\
        &= ((\pi^{0,1}\bar{\ot}{\rm id})\circ\nabla)\bar{\ot}{\rm id}+\sigma_{12}({\rm id}\bar{\ot}\pi^{0,1}\bar{\ot}{\rm id})(1\bar{\ot}{\nabla})\\
        &= \overline{\nabla} \bar{\ot} {\rm id} + \overline{\sigma}_{12} ({\rm id} \bar{\ot} \overline{\nabla}),
      \end{aligned}
    \end{eqnarray*}
    which is $\overline{\nabla}^{(\otimes2)}$.
\end{proof}
\begin{proposition}
    \label{descendsform}
    $\overline{\nabla}^{(\otimes2)}$ descends to a well defined left $\overline{\partial}$-operator on the bimodule $\Omega^{2}(\cla)$.
\end{proposition}
\begin{proof}
    We begin by recalling that for $\omega\in\Omega^{1}(\cla)$, $\nabla(\omega)=\theta\bar{\ot}\omega-\sigma(\omega\bar{\ot}\theta)$. Then for $\omega_{1},\omega_{2}\in\Omega^{1}(\cla)$, 
    \begin{eqnarray*}
      \begin{aligned}
        &\nabla^{(\otimes2)}(\omega_{1}\bar{\ot}\omega_{2})\\
&=\nabla(\omega_{1})\bar{\ot}\omega_{2}+\sigma_{12}(\omega_{1}\bar{\ot}\nabla(\omega_{2}))\\
        &= (\theta\bar{\ot}\omega_{1}\bar{\ot}\omega_{2})-\sigma_{12}(\omega_{1}\bar{\ot}\theta\bar{\ot}\omega_{2})+\sigma_{12}(\omega_{1}\bar{\ot}\theta\bar{\ot}\omega_{2})-\sigma_{12}\sigma_{23}(\omega_{1}\bar{\ot}\omega_{2}\bar{\ot}\theta)\\
        &= \theta\bar{\ot}\omega_{1}\bar{\ot}\omega_{2}-\sigma_{12}\sigma_{23}(\omega_{1}\bar{\ot}\omega_{2}\bar{\ot}\theta).
      \end{aligned}
    \end{eqnarray*}
    Now for $\omega_{1},\omega_{2}\in\Omega^{1}(\cla)$ such that $\sigma(\omega_{1}\bar{\ot}\omega_{2})=(\omega_{1}\bar{\ot}\omega_{2})$, using Lemma (\ref{relationnablanablabar}) and the braid relation of $\sigma$, we get
    \begin{eqnarray*}
      \begin{aligned}
        &\sigma_{23}\overline{\nabla}^{(\otimes2)}(\omega_{1}\bar{\ot}\omega_{2})\\
        &\quad=\sigma_{23}(\pi^{0,1}\bar{\ot}{\rm id}\bar{\ot}{\rm id})\nabla^{(\otimes2)}(\omega_{1}\bar{\ot}\omega_{2})\\
        &\quad= (\pi^{0,1}\bar{\ot}{\rm id}\bar{\ot}{\rm id})\sigma_{23}(\theta\bar{\ot}\omega_{1}\bar{\ot}\omega_{2}-\sigma_{12}\sigma_{23}(\omega_{1}\bar{\ot}\omega_{2}\bar{\ot}\theta))\\
        &\quad= (\pi^{0,1}\bar{\ot}{\rm id}\bar{\ot}{\rm id})(\theta\bar{\ot}\sigma(\omega_{1}\bar{\ot}\omega_{2})-\sigma_{23}\sigma_{12}\sigma_{23}(\omega_{1}\bar{\ot}\omega_{2}\bar{\ot}\theta))\\
        &\quad= (\pi^{0,1}\bar{\ot}{\rm id}\bar{\ot}{\rm id})(\theta\bar{\ot}\omega_{1}\bar{\ot}\omega_{2}-\sigma_{12}\sigma_{23}\sigma_{12}(\omega_{1}\bar{\ot}\omega_{2}\bar{\ot}\theta))\\
        &\quad= (\pi^{0,1}\bar{\ot}{\rm id}\bar{\ot}{\rm id})(\theta\bar{\ot}\omega_{1}\bar{\ot}\omega_{2}-\sigma_{12}\sigma_{23}(\omega_{1}\bar{\ot}\omega_{2}\bar{\ot}\theta)),
      \end{aligned}
    \end{eqnarray*}
    which is equal to $\overline{\nabla}^{(\otimes2)}(\omega_{1}\bar{\ot}\omega_{2})$. Therefore $\overline{\nabla}^{(\otimes2)}$ descends to a well defined $\mathbb{C}$-linear map from $\Omega^{2}(\cla)$ to $\Omega^{0,1}\bar{\ot}\Omega^{2}(\cla)$. We denote the descended map by $\overline{\nabla}^{(2)}$ which acts on an element $\omega_{1}\wedge\omega_{2}\in\Omega^{2}(\cla)$ by
    \begin{eqnarray}
    \label{nablawedge}
       \overline{\nabla}^{(2)}(\omega_{1}\wedge\omega_{2})=({\rm id}\bar{\ot}\wedge)\circ\overline{\nabla}^{(\otimes2)}(\omega_{1}\bar{\ot}\omega_{2}). 
    \end{eqnarray}
            To show that $\overline{\nabla}^{(2)}$ is again a left $\overline{\partial}$-operator, we use the facts that $({\rm id}\bar{\ot}\wedge)$ is a bimodule map, $\overline{\nabla}^{(\otimes 2)}$ is a left $\overline{\partial}$-operator on $\Omega^{1}(\cla)\bar{\ot}\Omega^{1}(\cla)$ and do the following calculation for $a\in\cla$,
    \begin{eqnarray*}
      \begin{aligned}
        &\overline{\nabla}^{(2)}(a.\omega_{1}\wedge\omega_{2})\\
        &\quad= ({\rm id}\bar{\ot}\wedge)\circ\overline{\nabla}^{(\otimes2)}(a.\omega_{1}\bar{\ot}\omega_{2})\\
        &\quad= ({\rm id}\bar{\ot}\wedge)(a\overline{\nabla}^{(\otimes2)}(\omega_{1}\bar{\ot}\omega_{2})+\overline{\partial}a\bar{\ot}\omega_{1}\bar{\ot}\omega_{2})\\
        &\quad= a({\rm id}\bar{\ot}\wedge)\circ\overline{\nabla}^{(\otimes2)}(\omega_{1}\bar{\ot}\omega_{2})+\overline{\partial}a\bar{\ot}\omega_{1}\wedge\omega_{2}\\
        &\quad= a\overline{\nabla}^{(2)}(\omega_{1}\wedge\omega_{2})+\overline{\partial}a\bar{\ot}\omega_{1}\wedge\omega_{2}.
      \end{aligned}
    \end{eqnarray*}
\end{proof}
\begin{remark}
    Note that in the above we haven't discussed the bimodule property of $\overline{\nabla}^{(2)}$. It can be shown to be a bimodule $\overline{\partial}$-operator with respect to a natural isomorphism $\psi:\Omega^{2}(\cla)\bar{\ot}\Omega^{0,1}(\cla)\raro\Omega^{0,1}(\cla)\bar{\ot}\Omega^{2}(\cla)$. But as we are interested in the holomorphic structure on the `exterior bundle' and $\Omega^{n}(\cla)=0$ for all $n\geq 3$, it suffices to consider $\overline{\nabla}^{(2)}$ as a left $\overline{\partial}$-operator.
\end{remark}
As observed earlier, with the almost complex structure $J$ on the $\ast$-differential calculus $(\Omega^{\bullet}(\cla),\der,\ast)$, $\Omega^{0,2}(\cla)=0$ and consequently for any left module $\cle$ over $\cla$, any left $\overline{\partial}$-operator $\overline{\nabla}$ induces a holomorphic structure canonically. Therefore $(\Omega^{1}(\cla),\overline{\nabla})$ and $(\Omega^{2}(\cla),\overline{\nabla}^{(2)})$ are holomorphic modules. Next we shall compute the holomorphic sections of the `bundles' $\Omega^{1}(\cla)$ and $\Omega^{2}(\cla)$.
\begin{proposition}
    \label{holosections}
    \begin{eqnarray*}
    && H^{0}(\Omega^{1}(\cla),\overline{\nabla})=\mathbb{C}\oplus\mathbb{C}\\
    && H^{0}(\Omega^{2}(\cla),\overline{\nabla}^{(2)})=\mathbb{C}
    \end{eqnarray*}
\end{proposition}
\begin{proof}
    Let $\omega=\sum_{\mu}(\alpha^{\mu}\xi_{\mu\raro\mu+1}+\beta^{\mu}\xi_{\mu\raro\mu-1})\in{\rm Ker}(\overline{\nabla})$ for some scalars $\{\alpha^{\mu},\beta^{\mu}\}_{\mu=1}^{n}$. Then by the formulae \eqref{nablabarformula}, we get
    \begin{equation*}
     \begin{aligned}
    \sum_{\mu}(\alpha^{\mu}(\xi_{\mu+1\raro\mu}&\bar{\ot}\xi_{\mu\raro\mu+1}-\xi_{\mu\raro\mu-1}\bar{\ot}\xi_{\mu-1\raro\mu})\\
    &+\beta^{\mu}(\xi_{\mu+1\raro\mu}\bar{\ot}\xi_{\mu\raro\mu-1}-\xi_{\mu\raro\mu-1}\bar{\ot}\xi_{\mu-1\raro\mu-2}))=0.
     \end{aligned}
    \end{equation*}
    Therefore comparing coefficients, we get $\alpha^{1}=\alpha^{2}=\cdots=\alpha^{n}$ and $\beta^{1}=\beta^{2}=\cdots=\beta^{n}$ proving that $H^{0}(\Omega^{1}(\cla),\overline{\nabla})={\rm Sp}_{\mathbb{C}}\{\sum_{\mu}\xi_{\mu\raro\mu+1},\sum_{\mu}\xi_{\mu\raro\mu-1}\}$.

    To compute $H^{0}(\Omega^{2}(\cla),\overline{\nabla}^{(2)})$, we write down the formulae for $\overline{\nabla}^{(2)}$ on the basis elements of $\Omega^{2}(\cla)$. Recall that a $\mathbb{C}$-linear basis for $\Omega^{2}(\cla)$ is given by $\{\xi_{\mu\raro\mu+1}\wedge\xi_{\mu+1\raro\mu}\}_{\mu=1,\ldots,n}$.
    \begin{eqnarray*}
      \begin{aligned}
        &\overline{\nabla}^{(2)}(\xi_{\mu\raro\mu+1}\wedge\xi_{\mu+1\raro\mu})\\
        &= ({\rm id}\bar{\ot}\wedge)\biggl(\overline{\nabla}(\xi_{\mu\raro\mu+1})\bar{\ot}\xi_{\mu+1\raro\mu}+\overline{\sigma}_{12}(\xi_{\mu\raro\mu+1}\bar{\ot}\overline{\nabla}(\xi_{\mu+1\raro\mu}))\biggr)\\
        &=({\rm id}\bar{\ot}\wedge)(\xi_{\mu+1\raro\mu}\bar{\ot}\xi_{\mu\raro\mu+1}\bar{\ot}\xi_{\mu+1\raro\mu}-\overline{\sigma}_{12}(\xi_{\mu\raro\mu+1}\bar{\ot}\xi_{\mu+1\raro\mu}\bar{\ot}\xi_{\mu\raro\mu-1}))\\
        &=({\rm id}\bar{\ot}\wedge)(\xi_{\mu+1\raro\mu}\bar{\ot}\xi_{\mu\raro\mu+1}\bar{\ot}\xi_{\mu+1\raro\mu}-\xi_{\mu\raro\mu-1}\bar{\ot}\xi_{\mu-1\raro\mu}\bar{\ot}\xi_{\mu\raro\mu-1})\\
        &=(\xi_{\mu+1\raro\mu}\bar{\ot}\xi_{\mu\raro\mu+1}\wedge\xi_{\mu+1\raro\mu}-\xi_{\mu\raro\mu-1}\bar{\ot}\xi_{\mu-1\raro\mu}\wedge\xi_{\mu\raro\mu-1}).
      \end{aligned}
    \end{eqnarray*}
    Therefore for an element $\omega=\sum_{\mu}\alpha^{\mu}\xi_{\mu\raro\mu+1}\wedge\xi_{\mu+1\raro\mu}$ such that $\omega\in{\rm Ker}(\overline{\nabla}^{(2)})$, we have
    \begin{displaymath}
        \sum_{\mu}\alpha^{\mu}(\xi_{\mu+1\raro\mu}\bar{\ot}\xi_{\mu\raro\mu+1}\wedge\xi_{\mu+1\raro\mu}-\xi_{\mu\raro\mu-1}\bar{\ot}\xi_{\mu-1\raro\mu}\wedge\xi_{\mu\raro\mu-1})=0.
    \end{displaymath}
    Again comparing coefficients like before, we obtain $\alpha^{1}=\alpha^{2}=\cdots=\alpha^{n}$ showing that $H^{0}(\Omega^{2}(\cla),\overline{\nabla}^{(2)})={\rm Sp}_{\mathbb{C}}\{\sum_{\mu}\xi_{\mu\raro\mu+1}\wedge\xi_{\mu+1\raro\mu}\}$ which proves the proposition.
\end{proof}
From Equation (\ref{nablawedge}), it is clear that if $\omega_{1},\omega_{2}\in{\rm Ker}(\overline{\nabla})$, then $\omega_{1}\wedge\omega_{2}\in{\rm Ker}(\overline{\nabla}^{(2)})$. As the holomorphic part of the polygon is connected, $H^{0}(\cla,\overline{\partial})={\rm Sp}_{\mathbb{C}}\{{\bf 1}\}$. Therefore $R:=H^{0}(\cla,\overline{\partial})\oplus H^{0}(\Omega^{1}(\cla),\overline{\nabla})\oplus H^{0}(\Omega^{2}(\cla),\overline{\nabla}^{(2)})$ is a unital ring with respect to the wedge product with ${\bf 1}\in H^{0}(\cla,\overline{\partial})$ being the unit of the ring.  We call $R$ the ring of `holomorphic sections of the exterior bundle'.
\begin{theorem}
    \label{Ring structute}
    The ring $R$ is isomorphic to the exterior algebra $\wedge^{\bullet}(\mathbb{C}^{2})$.
\end{theorem}
\begin{proof}
    Denoting $\sum_{\mu}\xi_{\mu\raro\mu+1}, \sum_{\mu}\xi_{\mu\raro\mu-1}$ by $X_{1},X_{2}$ respectively, we see that 
    \begin{displaymath}H^{0}(\Omega^{1}(\cla),\overline{\nabla})={\rm Sp}_{\mathbb{C}}\{X_{1},X_{2}\}.
    \end{displaymath}
    It is easy to verify that $X_{i}\wedge X_{i}=0$ for $i=1,2$ and $X_{1}\wedge X_{2}=-X_{2}\wedge X_{1}$. Moreover, $H^{0}(\Omega^{2}(\cla),\overline{\nabla}^{(2)})={\rm Sp}_{\mathbb{C}}\{X_{1}\wedge X_{2}\}$. Therefore denoting the standard basis of $\mathbb{C}^2$ by $e_{1},e_{2}$ and mapping $X_{i}$ to $e_{i}$ for $i=1,2$ and ${\bf 1}$ to $1\in\mathbb{C}$ establishes the isomorphism between $R$ and $\wedge^{\bullet}(\mathbb{C}^2)$. 
\end{proof}
\subsection{A positive Hochschild 2-cocycle and all that} 
In this section we shall show the existence of a positive 
Hochschild $2$-cocycle on the algebra of functions on $n$ points arising from 
the integrable almost complex structure $(\Omega^{\bullet}(\cla),\der,\ast,J)$ 
corresponding to the polygon graph as discussed in the previous subsection. 
The results of this subsection resembles those of \cite{qbundle} although we 
do not consider any twisted versions of cyclic or Hochschild cohomology. 
A few quick definitions are in order. For generalities on Hochschild and cyclic 
cohomology the reader is referred to Section 3.1 of chapter 3 of \cite{Khalkhali}.
\begin{definition}
    Let $(C^{\bullet}(\cla),b)$ be the Hochschild complex of $\cla$ over the field of complex numbers i.e.
    \begin{displaymath}
     C^{n}(\cla):={\rm Hom}(\cla^{\otimes(n+1)},\mathbb{C})\end{displaymath}
      \begin{eqnarray*}
      b\phi(a_{0},a_{1},\ldots,a_{n+1})&:=& \sum_{j=0}^{n}(-1)^{j}\phi(a_{0},\ldots,a_{j}a_{j+1},\ldots,a_{n+1})\\
     &&\quad+(-1)^{n+1}\phi(a_{n+1}a_{0},a_{1},\ldots,a_{n}).
    \end{eqnarray*}
  An element $\phi\in C^{n}(\cla)$ is called a Hochschild $n$-cochain. An $n$-cochain $\phi$ is said to be cyclic if $\phi(a_{n},a_{0},\ldots,a_{n-1})=(-1)^{n}\phi(a_{0},\ldots,a_{n})$.      
\end{definition}
\begin{definition}
    A Hochschild $2n$-cocycle on a $\ast$-algebra $\cla$ is said to be positive if the following pairing defines a positive sesquilinear form on $\cla^{\otimes(n+1)}$:
    \begin{displaymath}
        \langle a_{0}\ot\cdots\ot a_{n}, b_{0}\ot\cdots\ot b_{n}\rangle=\phi(b_{0}^{\ast}a_{0},a_{1},\ldots,a_{n},b_{n}^{\ast},\ldots,b_{1}^{\ast})
    \end{displaymath}
\end{definition}

Recall the closed graded trace $\int$ of dimension $2$ on 
the $2$-dimensional $\ast$-differential calculus $(\Omega^{\bullet}(\cla),\der,\ast)$ from Lemma \ref{gradedtrace} where $\cla$ is the $\ast$-algebra of functions on $n$-points. 
\begin{lemma}
    \label{Hcocylce}
    The following formula produces a positive Hochschild $2$-cocycle on $\cla$:
    \begin{displaymath}
        \phi(f_{0},f_{1},f_{2}):=\int f_{0}\partial(f_{1})\wedge\overline{\partial}(f_{2}).
    \end{displaymath}
\end{lemma}
\begin{proof}
    The equality $\int f_{0}\partial(f_{1})\wedge\overline{\partial}(f_{2})f_{3}=\int f_{3}f_{0}\partial(f_{1})\wedge\overline{\partial}(f_{2})$ follows easily from the previously observed fact that the right and left action of a function on the vertex set on the bimodule $\Omega^{2}(\cla)$ coincide.
    Then using the Leibniz rule for $\partial$ and $\overline{\partial}$, $\phi$ is shown to be a $2$-cocycle along the lines of~\cite{qbundle}*{Proposition 5.2}. To show the positivity, we compute $\phi(f_{0}^{\ast}f_{0},f_{1},f_{1}^{\ast})$ for $f_{0},f_{1}\in\cla$.
    \begin{eqnarray*}
      \begin{aligned}
        &\int f_{0}^{\ast}f_{0}\partial(f_{1})\wedge\overline{\partial}(f_{1}^{\ast})\\
        &=\int f_{0}  \partial(f_{1})\wedge (f_{0}\partial(f_{1}))^{\ast}\\
        &=- \int \sum_{\mu=1}^{n}f_{0}(\mu)(f_{1}(\mu+1)-f_{1}(\mu))\xi_{\mu\raro\mu+1}
        \wedge\sum_{\nu=1}^{n} \overline{f_{0}(\nu)}\overline{(f_{1}(\nu+1)-f_{1}(\mu))}\xi_{\nu+1\raro\nu}\\ 
        &=-\int \sum_{\mu=1}^{n}|f_{0}(\mu)|^{2}|f_{1}(\mu+1)-f_{1}(\mu)|^{2}\xi_{\mu\raro\mu+1}\wedge\xi_{\mu+1\raro\mu}\\
        &=\int \sum_{\mu=1}^{n}|f_{0}(\mu)|^{2}|f_{1}(\mu+1)-f_{1}(\mu)|^{2}\xi_{\mu\raro\mu-1}\wedge\xi_{\mu-1\raro\mu}\\
        &= \frac{1}{n}\sum_{\mu=1}^{n}|f_{0}(\mu)|^{2}|f_{1}(\mu+1)-f_{1}(\mu)|^{2}\geq 0.
      \end{aligned}
    \end{eqnarray*}
\end{proof}
Recall that the closed graded trace $\int$ of dimension $2$ gives rise to the cyclic $2$-cocycle on $C(V)$ given by the formula (see \cite{Khalkhali})
\begin{displaymath}
    \tau(f_{0},f_{1},f_{2})=\frac{1}{2}\int f_{0}\der{f_{1}}\wedge \der{f_{2}}.
\end{displaymath}
Unfortunately this cocycle turns out to be a trivial one. We omit the proof of triviality of $\tau$. One can check directly by observing the action of $\tau$ on the basis elements of $C(V\times V\times V)$ which is a finite dimensional vector space. For example, for any vertex $\mu$, $\tau(\delta_{\mu}\ot\delta_{\mu}\ot\delta_{\mu})=0$. The non trivial actions of $\tau$ are on the basis elements of type $\delta_{\mu}\ot\delta_{\mu-1}\ot\delta_{\mu}$. But then $\tau=b\tau^{\prime}$ for a cyclic $1$-cycle $\tau^{\prime}$ that takes the value $1$ on the basis elements $\{\delta_{\mu}\ot\delta_{\mu+1}\}_{\mu=1}^{n}$ which renders $\tau$ to be trivial.
\begin{proposition}
    The cocycles $\tau$ and $\phi$ are cohomologous and consequently the positive Hochschild $2$-cocycle $\phi$ is trivial.
\end{proposition}
\begin{proof}
    We shall adopt the proof of~\cite{qbundle}*{Proposition 5.4}. So we define a Hochschild $1$-cochain $\psi$ on $C(V)$ by the following:
    \[
        \psi(f_{0},f_{1})=\frac{1}{2}\int f_{0}\partial\overline{\partial}(f_{1}).
    \]
    Then following the same calculation as in~\cite{qbundle}*{Proposition 5.4}, we have
    \begin{displaymath}
        b\psi(f_{0},f_{1},f_{2})=\frac{1}{2}\int\biggl(f_{0}f_{1}\partial\overline{\partial}(f_{2})-f_{0}\partial\overline{\partial}(f_{1}f_{2})+f_{0}\partial\overline{\partial}(f_{1})f_{2}\biggr).
    \end{displaymath}
    On the other hand, using the fact that $\partial(f_{1})\wedge\partial(f_{2})=\overline{\partial}(f_{1})\wedge\overline{\partial}(f_{2})=0$ we obtain
    \begin{displaymath}
        \frac{1}{2}\int f_{0}\der{f_{1}}\wedge \der{f_{2}}-\int f_{0}\partial f_{1}\wedge\overline{\partial}f_{2}=\frac{1}{2}\int\bigl(f_{0}\overline{\partial}f_{1}\wedge\partial f_{2}-f_{0}\partial f_{1}\wedge\overline{\partial}f_{2}\bigr).
    \end{displaymath}
    Now using the graded Leibniz rule for $\partial$ (see \cite{Beggs}), we get
    \begin{eqnarray*}
        -\partial\overline{\partial}(f_{1}f_{2})+f_{1}\partial\overline{\partial}(f_{2})+\partial\overline{\partial}(f_{1})f_{2}=\bigl(f_{0}\overline{\partial}f_{1}\wedge\partial f_{2}-f_{0}\partial f_{1}\wedge\overline{\partial}f_{2}\bigr),
    \end{eqnarray*}
    showing that
    \begin{displaymath}
       \frac{1}{2}\int f_{0}\der{f_{1}}\wedge \der{f_{2}}=\int f_{0}\partial f_{1}\wedge\overline{\partial}f_{2} + \frac{1}{2}\int f_{0}\bigl(-\partial\overline{\partial}(f_{1}f_{2})+f_{1}\partial\overline{\partial}(f_{2})+\partial\overline{\partial}(f_{1})f_{2}\bigr). 
    \end{displaymath}
Therefore, $(\tau-\phi)(f_{0},f_{1},f_{2})=b\psi(f_{0},f_{1},f_{2})$ proving that $\tau$ and $\phi$ are cohomologues cocycles. In particular $\phi$ is a trivial Hochschild $2$-cocycle.
\end{proof}

\section*{Concluding remarks}\begin{itemize}
\item In the particular example of polygon, we have studied a holomorphic structure induced by a bimodule connection on the module of one-forms. As observed in \cite{Majid}, such bimodule connections are in on-one correspondence with a pair of bimodule maps $(\sigma,\alpha)$ where $\sigma:\Omega^{1}(\cla)\bar{\ot}\Omega^{1}(\cla)\raro \Omega^{1}(\cla)\bar{\ot}\Omega^{1}(\cla)$ and $\alpha:\Omega^{1}(\cla)\raro \Omega^{1}(\cla)\bar{\ot}\Omega^{1}(\cla)$. We have studied the bimodule connection with respect to $(\sigma,0)$ where $\sigma$ is the bimodule map giving the canonical prolongation. Keeping $\sigma$ fixed and varying $\alpha$, it will be interesting to study equivalence of different holomorphic structures coming from different $\alpha$.
\item We have not discussed about any K$\ddot{{\rm a}}$hler structure (see \cites{kahlerbuachala, Das}) in this paper. We hope to investigate K$\ddot{{\rm a}}$hler structure and consequently relevant Hodge theory on bidirected graphs in near future.\end{itemize}

\end{document}